\documentclass{amsart}[1996/10/24] 

\newlength\skwd\newcommand\sk[1]{\hskip#1\skwd}
\newcommand\mymail{\settowidth\skwd{x}%
\mbox{\sk{5}@\sk{2}b\sk{-9}d\sk{9}d\sk{-9}e\sk{4}e%
\sk{3}e\sk{-9}l\sk{-3}w\sk{4}w\sk{-3}z\sk{4}.\sk{2}}}

\newcommand\xxv{\settowidth\skwd{9}%
\mbox{\sk{1}5\sk{-2}2\sk{1}}}
\newcommand\myzip{\settowidth\skwd{9}%
\mbox{\sk{3}02\sk{-5}588\sk{2}}}

\DeclareFontFamily{T1}{wncysc}{\hyphenchar\font=-1}
\DeclareFontShape{T1}{wncysc}{m}{sc}{<->wncysc10}{} 
\newcommand\cysc{\usefont{T1}{wncysc}{m}{sc}}

\newcommand{\cyscB}{{\cysc{V}}} 
\newcommand{\cyscb}{{\cysc{v}}} 

\newcommand\sss{\sffamily\slshape}

\newcommand\ic{\/\strut}

\newcommand\hh{\nobreak-\allowbreak\hskip0pt\relax}
\newcommand\hd{\nobreak--\allowbreak\hskip0pt\relax}
\newcommand\hs{\nobreak/\allowbreak\hskip0pt\relax}

\newcommand\singlenumber[2]{\expandafter\xdef\csname the#1\endcsname
{\expandafter\noexpand\csname the#2\endcsname}}

\let\amsshoveleft\shoveleft
\renewcommand\shoveleft[1]{\amsshoveleft{#1\hfill}}

\DeclareMathOperator{\artanh}{artanh} 

\newcommand\np{\!+\!} 

\DeclareFontFamily{T1}{cmtex}{\hyphenchar\font=-1}
\DeclareFontShape{T1}{cmtex}{m}{n}{<-9>cmtex8<9>cmtex9<10->cmtex10}{}
\newcommand\mtt{\usefont{T1}{cmtex}{m}{n}}

\newcommand{\tneq}{{\mtt\char26}} 

\newcommand\xtra{\catcode`^12\catcode`_12}

{\obeyspaces\gdef\sqze{\obeyspaces\let \thinspace}}

\newcommand\ttt{\tt\xtra\sqze} 

\newcommand\tttdef[1]{%
  \long\def\gobble ##1{}
  \edef\tttdefa{\expandafter\gobble\string#1}%
  \edef#1{\noexpand\protect 
    \expandafter\noexpand\csname\tttdefa\space\endcsname}%
  \begingroup
    \def\tttdefb##1{%
      \expandafter\xdef\csname\tttdefa\space\endcsname{%
        {\noexpand\tt##1}}}
    \xtra\sqze 
    \afterassignment\endgroup
    \tttdefb}

\newcommand\fb{\allowbreak\hskip0pt\relax} 

\newcommand\ads{\vskip\abovedisplayskip\noindent}
\newcommand\bds{\vskip\belowdisplayskip\noindent}

\newlength{\lind}

\newcommand\reu[1]{\hspace*\multlinegap\hskip#1\lind\ttt}
\newcommand\len[1]{\refstepcounter{equation}(\theequation)
  \hskip\multlinetaggap\hskip#1\lind\ttt}
\newcommand\leu[1]{\phantom{(\theequation)}
  \hskip\multlinetaggap\hskip#1\lind\ttt}

\newcommand\bfl{\begingroup\raggedright}
\newcommand\efl{\par\endgroup} 

\begin{document}

\title[Recurrence formulae for algebraic integrals]
{Two-term recurrence formulae for\\ indefinite algebraic integrals}


\author{Detmar Martin Welz}
\address{Unterm E{\cyscb}erg {\xxv}, {\myzip} {\cyscB}alve, Germany} 
\email{\mymail} 


\date{\today}



\begin{abstract}
Two-term recurrence relations are supplied for indefinite integrals of 
functions that involve factors of the types ${P_2}^n$, ${P_3}^n$, 
${P_4}^n$, ${P_1}^m\,{Q_1}^n$, $E_1 {P_1}^n$, ${P_1}^m\,{Q_2}^n$, 
$E_1 {P_2}^n$, ${P_2}^m\,{Q_2}^n$, ${P_1}^m\,{Q_1}^n\,{S_1}^p$, 
$E_1 {P_1}^m\,{Q_1}^n$, ${P_1}^m\,{Q_1}^n\,{S_2}^p$, and 
${P_1}^m\,{Q_1}^n\,{S_1}^p\,{T_1}^q$, where $P_i$, $Q_j$, $S_k$ and $T_l$ 
denote arbitrary polynomials of degree $i$, $j$, $k$ and $l$ in the 
integration variable, $E_1$ represents the exponential function of an 
arbitrary linear polynomial in this variable, and $m$, $n$, $p$ and $q$ 
are arbitrary constant exponents. The $136$ relations leave the form of 
an integrand unchanged and increment or decrement the exponents in steps 
of unity.
\end{abstract}

\maketitle

\section*{Introduction}

To be maximally useful in science and engineering applications, the 
symbolic evaluation of an indefinite integral should be expressed as 
compactly as possible and in a form suitable for verification by 
differentiation, and hold in the entire complex plane whenever feasible. 
With free parameters in the integrand taken to be real, it should further 
be continuous on the real line where the integrand is integrable along 
it, and preferably be real on this line if the integrand is real 
everywhere on it, and also not explicitly involve the imaginary unit 
unless the integrand does too. When these goals are in partial conflict, 
the resolution should be sought by sacrificing compactness first, absence 
of imaginary offset and imaginary unit next, and continuity at integrable 
poles last.

These principles will find, if not adoption, then ready approval. Indeed, 
the vast majority of the integrands arising in practice do not involve 
the imaginary unit or complex parameters (or else they easily separate 
into real and imaginary parts) and their antiderivatives are needed for 
real values of the integration variable only, where the absence of 
unexpected, symbolically often hard to localize discontinuities ensures 
that definite integrals calculated according to Newton and Leibniz will 
be correct. By default, an integrator should therefore place the 
unavoidable branch cuts such in an antiderivative that they do not 
traverse the real line, employing to this end functional equations for 
the logarithmic and related functions implicated in those cases. Note 
that the details may depend on the actual function implementation.

To the recurrent chagrin of countless users, the automated integrators of 
today's computer{\hh}algebra systems do not systematically respect these 
principles. Consider the parameterless example of
\begin{multline*}
\int \!\frac {1 - x + 3 x^2} {(1 + x + x^2)^2 \,\sqrt{1 - x + x^2}} 
  \;\mathrm{d}x\\
\shoveleft{ = \frac {(1 + x) \sqrt{1 - x + x^2}} {1 + x + x^2} 
  + \frac{1}{2} \int \!\frac {3 - x} {(1 + x + x^2) \sqrt{1 - x + x^2}} 
  \;\mathrm{d}x }\\
\shoveleft{ = \frac {(1 + x) \sqrt{1 - x + x^2}} {1 + x + x^2}
  + \sqrt{2} \arctan \frac {\sqrt{2} \,(1 + x)} {\sqrt{1 - x + x^2}} 
  - \frac{1}{\sqrt{6}} \artanh \frac {\sqrt{6} \,(1 - x)} 
  {3 \sqrt{1 - x + x^2}}, }
\end{multline*}
in which the integrand is analytic everywhere on the real line. The first 
reduction step invokes a recurrence formula to raise the exponent $-2$ in 
the integrand to $-1$, while the second step employs transformations and 
evaluations chosen for compliance with the principles stated. The major 
commercial contenders, {\sc{Maple}} (version 16) and {\sc{Mathematica}} 
(version 8), here produce discontinuous antiderivatives of about ten 
times this size; {\sc{Maple}} is able to correct for the discontinuity 
(an imaginary step at $x = 1$), whereas {\sc{Mathematica}} bungles 
definite integrals that straddle the discontinuity (a real jump near 
$x = -0.493$).

The short example already illustrates the potential of combining suitable 
recurrence relations with carefully chosen terminal evaluations to 
achieve fast and strong automated symbolic integration in accordance with 
the above principles. This approach is currently being explored in the 
rule-based integrator {\sc{Rubi}} \cite{JR}, for instance. The purpose of 
the present publication is to furnish efficient recurrence formulae for 
the reduction of algebraic integrals that are regularly encountered in 
practice; recipes for a properly principled evaluation of the terminal 
instances of the integrals are to be published separately.

Being two-term relations, the recurrence formulae presented here replace 
an algebraic integrand with a simpler one while liberating a 
corresponding algebraic part of the antiderivative. Specifically, they 
apply to integrands that involve products of arbitrary powers of linear 
or quadratic polynomials, or isolated arbitrary powers of cubic or 
quartic polynomials; this naturally includes integer positive and 
negative powers. The relations increment or decrement the exponents in 
steps of unity, either separately or in up-down pairs, and leave the form 
of an integrand unchanged. Accordingly, they should lend themselves 
especially to the automated symbolic integration of algebraic functions 
in which polynomials are raised to high positive or negative powers.

Also covered are integrands containing the exponential factor 
$\exp(a + b x)$ along with such powers of polynomials; this exponential 
derives from the power of a linear polynomial via Euler's famous limit: 
$\exp(a + b x) = \lim_{n\to\infty}(1 + a/n + b/n\:x)^n$. Insofar as 
simple analog formulae for the trigonometric factors $\cos(a + b x)$ and 
$\sin(a + b x)$ exist, they are supplied as well; otherwise these 
integrands should be split such that the integrals involve either 
$\exp[\,i \,(a + b x)]$ or $\exp[-i \,(a + b x)]$, where $\pm i a$ and 
$\pm i b$ replace the parameters $a$ and $b$. Note that the exponential 
factors cannot be recurred on because their parameters are unrelated to 
the exponent of the linear polynomial.

Quite generally, indefinite integrals of the present kind with generic 
exponents evaluate to special cases of those hypergeometric functions 
that admit one{\hh}dimensional Euler-type integral representations. With 
the sum of polynomial degrees in the algebraic integrands ranging from 
two to four, the Gauss function $_2F_1$, the Appell function $F_1$, and 
the Lauricella function $F_D$ of three variables result, and where an 
exponential factor accompanies polynomials with a combined degree of one 
or two, the incomplete gamma function (which is a special case of 
Kummer's confluent hypergeometric function) and the Humbert function 
$\varPhi_1$ appear. The recurrence relations for the integrals are thus 
equivalent to relations among contiguous instances of these 
hypergeometric functions, and they are of interest mainly for integrals 
whose evaluation reduces to elementary functions or to the canonical 
elliptic integrals.

Recurrence formulae for integrands that include the power of a linear 
polynomial, say $(a + b x)^m$, admit profitable transformations once this 
factor is specialized to the monomial $x^m$ by setting $a = 0$ and 
$b = 1$. Relation \eqref{9A1} from Section \ref{9A} below, for example, 
can thus be specialized to
\begin{multline*}
(n \np 1) (a d - b c) \,a \!\int (A + B x) 
  \,x^m (a + b x)^n (c + d x)^p \;\mathrm{d}x\\
\shoveleft{ {}- \!\int \!\big[(n \np 1) \,A\, (a d - b c) 
  - (m \np 1) (A b - B a) \,c 
  - (m \np n \np p \np 3) (A b - B a) \,d x\big] }\\[-1ex]
\shoveright{ x^m (a + b x)^{n + 1} (c + d x)^p \;\mathrm{d}x }\\
\shoveleft{ {}- (A b - B a) 
  \,x^{m + 1} (a + b x)^{n + 1} (c + d x)^{p + 1} = 0. }
\end{multline*}
The substitution $x \gets x^q$, $m \gets (m \np 1)/q - 1$ now yields
\begin{multline*}
(n \np 1) (a d - b c) \,a \!\int (A + B x^q) 
  \,x^m (a + b x^q)^n (c + d x^q)^p \;\mathrm{d}x\\
\shoveleft{ {}- \!\int \!\big[(n \np 1) \,A\, (a d - b c) 
  - \lambda\, (A b - B a) \,c 
  - (\lambda \np n \np p \np 2) (A b - B a) \,d x^q\big] }\\[-1ex]
\shoveright{ x^m (a + b x^q)^{n + 1} (c + d x^q)^p \;\mathrm{d}x }\\
\shoveleft{ {}- \frac{1}{q}\, (A b - B a) 
  \,x^{m + 1} (a + b x^q)^{n + 1} (c + d x^q)^{p + 1} = 0, }
\end{multline*}
where $\lambda = (m \np 1)/q$.
In the new kind of integrand, the monomial exponent can be changed only 
in steps of size $q$, and in the present relation the monomial can be 
eliminated altogether from the integrands by setting $m = 0$. Recurrence 
formulae of this type are needed for the reduction of elliptic integrals 
in particular, for which $q = 2$.

In a further step, the monomial power may be split by substituting $m 
\gets m + r n + s p$, and then $x^{r n}$ and $x^{s p}$ be merged with the 
respective other factors:
\begin{multline*}
(n \np 1) (a d - b c) \,a \!\int (A + B x^q) 
  \,x^m (a x^r + b x^{r + q})^n (c x^s + d x^{s + q})^p \;\mathrm{d}x \\
\shoveleft{ {}- \!\int \!\big[(n \np 1) \,A\, (a d - b c) 
  - \mu\, (A b - B a) \,c 
  - (\mu \np n \np p \np 2) (A b - B a) \,d x^q\big] }\\[-1ex]
\shoveright{ x^{m - r} (a x^r + b x^{r + q})^{n + 1} 
  (c x^s + d x^{s + q})^p \;\mathrm{d}x }\\
\shoveleft{ {}- \frac{1}{q}\, (A b - B a) 
  \,x^{m - r - s + 1} (a x^r + b x^{r + q})^{n + 1} 
  (c x^s + d x^{s + q})^{p + 1} = 0, }
\end{multline*}
where $\mu = (m + r n + s p + 1)/q$. The monomial exponent has now 
clearly lost its independence. Note that recurrence relations transformed 
in these ways hold for integrands involving arbitrary fractional 
exponents because possible branch-cut crossings affect all terms of a 
relation equally; thus $m$, $n$, $p$, $q$, $r$ and $s$ in the above are 
unrestricted apart from the proviso that $q \neq 0$. The same{\pagebreak[1]} 
transformations can obviously be applied to relations in which the 
integrand includes an exponential factor. The conversion of the formulae 
presented below, however, is left to the fortunate reader since this is a 
fairly straightforward task whose combined results would have consumed an 
unwarranted amount of space here.

The present formula collection is kept at a comfortable size by limiting 
the sum of polynomial degrees in the algebraic integrands at four; to 
just include, in the notation of the Abstract, the degree-five product 
${P_1}^m\,{Q_4}^n$ in the same way would roughly double the size already. 
Among the degree-four integrands, no combination is omitted but 
${P_1}^m\,{Q_3}^n$, as it could aid a systematic handling of elliptic 
integrals (for integer $m$ and odd integer $2 n$) alongside its cousin 
${P_1}^m\,{Q_4}^n$ only. The combined degree of powers that accompany an 
exponential (or trigonometric) factor of linear argument is even limited 
at two: the degree-three products $E_1 {P_3}^n$, $E_1 {P_1}^m\,{Q_2}^n$, 
and $E_1 {P_1}^m\,{Q_1}^n\,{S_1}^p$ already appear to be of no interest 
in practice, and their integrals do not relate to one of the established 
hypergeometric functions either. Exponential factors with higher-degree 
arguments are ignored for the same reason.

Still, when adapted relations are applied to suitably transformed 
integrals, the present formulae suffice to handle any elliptic integrand. 
Relations for the integrand factors $x^m\,{Q_2}^n$ and 
$x^m\,{Q_1}^n\,{S_2}^p$ with the substitution $x \gets x^2$ in particular 
encompass two-term equivalents of the multiterm recurrence formulae 
traditionally employed to reduce general elliptic integrals after the 
$x^1$ and $x^3$ terms of the radicand have been annihilated by a M\"obius 
transformation and the even part of the rational cofactor has been 
expanded into partial fractions with respect to $x^2$. Note that the 
M\"obius transformation generates an ugly piecewise{\hh}constant 
prefactor and that the expansion step may introduce complex coefficients, 
but these effects beset any reduction to canonical elliptic integrals.

Recurrence relations for rational integrands multiplied by a power 
${P_3}^n$ or ${P_4}^n$ are actually made redundant for any integer $4 n$ 
in this way; the direct reduction being much easier and more general, 
however, it is supported here for the isolated powers, although the 
terminal evaluation may still require a M\"obius transformation. 
Moreover, a piecewise{\hh}constant prefactor would arise for odd $2 n$ 
even if an isolated radicand has repeated roots, but is then better 
avoided, as it is rarely a natural part of the elementary antiderivative 
that results. The other degree-three and degree-four combinations of 
powers covered below must be regarded as indispensable anyway; so the 
elliptic cases among them are best reduced like the elementary ones, and 
the fractional powers just collected into a single factor prior to the 
terminal evaluation, introducing a piecewise constant which fuses with 
that from a subsequent M\"obius transformation.

For each integrand, some of the recurrence relations will fail if certain 
polynomial roots coincide, because the discriminants or resultants appear 
as factors in front of the integrals. Such confluent roots can be fused, 
and thereby made harmless, through a rational refactorization of the 
integrand, which lowers the total degree, but introduces ugly 
piecewise{\hh}constant prefactors in general. Dedicated formulae are 
therefore provided here for all cases of degeneracy from confluent 
polynomial roots, whereby the need for refactorization is avoided, or at 
least deferred until the terminal integral is evaluated. In formulating 
mutually exclusive conditions for the various cases of degeneracy below, 
the polynomial factors have been assumed to possess nonzero leading 
coefficients. The conditions and associated recurrence formulae are most 
easily verified on polynomials having symbolic roots of the various 
possible confluence patterns.

The two-term recurrence relations have been derived by the method of 
undetermined coefficients and can be checked simply by differentiation. 
To facilitate their use in computer{\hh}algebra systems, they are listed 
here in a linear functional notation without subscripts or superscripts, 
which should be self{\hh}explaining. Extensive abbreviation of repeated 
subexpressions on both the integrand and the formula level serves to 
condense the listing and to expose the structure of the relations. While 
two-term formulae for the simplemost integrands $(a + b x + c x^2)^n$, 
$(a + b x)^m (c + d x)^n$ and $\exp(a + b x) (c + d x)^n$ can be found in 
standard tables of integrals (so in Gradshteyn{\hd}Ryzhik \cite{GR} as 
item 2.260.2; as items 2.151, 2.153, 2.155; and for $a = c = 0$, $d = 1$ 
as items 2.321.1, 2.324.1), they are included below for convenience.

The state of the art prior to the present publication is reflected in a 
1948 book by Timofeyev \cite{Ti} (there is no evidence for the 1933 
edition cited in \cite{GR}) and a recent communication by Barnett 
\cite{Ba}.

\begingroup\interlinepenalty200\relax 
\singlenumber{equation}{section}
\tttdef\formi{(a + b*x + c*x^2)^n}
\section{Integrands involving \formi}

To exclude integrands with confluent roots, the following recurrence 
should be applied only if the quadratic discriminant does not vanish: 
{\ttt{4*a*c - b^2 \tneq 0}}. Note that the relation can be used backwards 
for its own inverse.
\ads
{\len0{(n + 1)*(4*a*c - b^2)*INT((a + b*x + c*x^2)^n, x)}}\\
{\leu0{ - 2*(2*n + 3)*c*INT((a + b*x + c*x^2)^(n + 1), x)}}\\
{\leu0{ + (b + 2*c*x)*(a + b*x + c*x^2)^(n + 1) = 0}}.

\tttdef\formii{(a + b*x + c*x^2 + d*x^3)^n}
\section{Integrands involving \formii}
\numberwithin{equation}{subsection}

Abbreviations used with these integrands:
\ads
{\reu0{ra = 2*(3*b*d - c^2)}}, {\ttt{rb = 9*a*d - b*c}}, 
{\ttt{rc = 2*(3*a*c - b^2)}}.

\subsection{Nondegenerate integrands}
\label{2A}

To exclude integrands with confluent roots, the following recurrences 
should be applied only if the cubic discriminant does not vanish: 
{\ttt{ra*rc - rb^2 \tneq 0}}.
\ads
{\len0{(n + 1)*(ra*rc - rb^2)*INT((A + B*x)*(a + b*x + c*x^2 + d*x^3)^n, x)}}\\
{\leu0{ + 3*INT(((3*n + 4)*U*d - (2*n + 3)*V*c - (3*n + 5)*V*d*x)}}\\
{\leu1{*(a + b*x + c*x^2 + d*x^3)^(n + 1), x)}}\\
{\leu0{ - (U*c - 2*V*b + 3*(U*d - V*c)*x - 3*V*d*x^2)}}\\
{\leu1{*(a + b*x + c*x^2 + d*x^3)^(n + 1) = 0}},
\bds\bfl
where {\ttt{U = A*rb - B*rc}}, {\ttt{V = A*ra - B*rb}}.
\efl\ads
{\len0{3*v*(v + d)*INT((A + B*x)*(a + b*x + c*x^2 + d*x^3)^n, x)}}\\
{\leu0{ - n*INT((U*rb + B*v*rc + (U*ra + B*v*rb)*x)}}\\
{\leu1{*(a + b*x + c*x^2 + d*x^3)^(n - 1), x)}}\\
{\leu0{ - (U*c + 2*B*v*b + 3*(U*d + B*v*c)*x + 3*B*v*d*x^2)}}\\
{\leu1{*(a + b*x + c*x^2 + d*x^3)^n = 0}},
\bds\bfl
where {\ttt{U = (3*n + 2)*A*d - (2*n + 1)*B*c}}, {\ttt{v = (3*n + 1)*d}}.
\efl

\singlenumber{equation}{subsection}
\subsection{Singly degenerate integrands}
\label{2B}

Here the cubic polynomial is required to possess a double root whence its{\pagebreak[3]} 
discriminant must vanish: {\ttt{ra*rc - rb^2 = 0}}. To exclude doubly 
degenerate cubics, the recurrence should be applied only if additionally 
{\ttt{ra \tneq 0}}. Note that the relation can be used backwards for its 
own inverse.
\ads
{\len0{(n + 1)*(2*n + 1)*(ra*c - 3*rb*d)*INT((a + b*x + c*x^2 + d*x^3)^n, x)}}\\
{\leu0{ + 3*(3*n + 4)*w*d*INT((a + b*x + c*x^2 + d*x^3)^(n + 1), x)}}\\
{\leu0{ - (n + 1)*(rb*c - 2*rc*d + (ra*c - rb*d)*x + ra*d*x^2)}}\\
{\leu1{*(a + b*x + c*x^2 + d*x^3)^n}}\\
{\leu0{ - w*(c + 3*d*x)*(a + b*x + c*x^2 + d*x^3)^(n + 1) = 0}},
\bds\bfl
where {\ttt{w = (3*n + 2)*d}}.
\efl

\tttdef\formiii{(a + b*x + c*x^2 + d*x^3 + e*x^4)^n}
\section{Integrands involving \formiii}
\numberwithin{equation}{subsection}

Abbreviations used with these integrands:
\ads
{\reu0{ra = 8*c*e - 3*d^2}}, {\ttt{rb = 6*b*e - c*d}}, 
{\ttt{rc1 = 16*a*e - b*d}},\\
{\reu0{rc2 = 4*a*e + 2*b*d - c^2}}, {\ttt{rd = 6*a*d - b*c}}, 
{\ttt{re = 8*a*c - 3*b^2}}.\\[1ex]
{\reu0{sa = ra*rc2 - rb^2}}, {\ttt{sb = 1/2 *(ra*rd - rb*rc1)}}, 
{\ttt{sc = rb*rd - rc1*rc2}},\\
{\reu0{sd = 1/2 *(rb*re - rc1*rd)}}, {\ttt{se = rc2*re - rd^2}}.

\subsection{Nondegenerate integrands}
\label{3A}

To exclude integrands with confluent roots, the following recurrences 
should be applied only if the quartic discriminant does not vanish: 
{\ttt{sa*re - 2*sb*rd + sc*rc1 \tneq 0}}.
\ads
{\len0{(n + 1)*(sa*re - 2*sb*rd + sc*rc1)*INT((A + B*x + C*x^2)}}\\
{\leu1{*(a + b*x + c*x^2 + d*x^3 + e*x^4)^n, x)}}\\
{\leu0{ - 4*INT(((4*n + 5)*U*e - (3*n + 4)*V*d + u*W*c}}\\
{\leu1{ - (2*u*V*e - (3*n + 5)*W*d)*x + (4*n + 7)*W*e*x^2)}}\\
{\leu1{*(a + b*x + c*x^2 + d*x^3 + e*x^4)^(n + 1), x)}}\\
{\leu0{ + (U*d - 2*V*c + 3*W*b + 4*(U*e - V*d + W*c)*x}}\\
{\leu1{ - 4*(V*e - W*d)*x^2 + 4*W*e*x^3)}}\\
{\leu1{*(a + b*x + c*x^2 + d*x^3 + e*x^4)^(n + 1) = 0}},
\bds\bfl
where {\ttt{U = A*sc - B*sd + C*se}}, {\ttt{V = A*sb - B*sc + C*sd}}, 
{\ttt{W = A*sa - B*sb + C*sc}}, {\ttt{u = 2*n + 3}}.
\efl\ads
{\len0{4*(2*n + 1)*(4*n + 3)*w*e^2*INT((A + B*x + C*x^2)}}\\
{\leu1{*(a + b*x + c*x^2 + d*x^3 + e*x^4)^n, x)}}\\
{\leu0{ - n*INT(((U - V*v)*rc1 + (V*rd + W*re)*w + 2*((U - V*v)*rb}}\\
{\leu1{ + (V*rc2 + W*rd)*w)*x + ((U - V*v)*ra + (V*rb + W*rc1)*w)*x^2)}}\\
{\leu1{*(a + b*x + c*x^2 + d*x^3 + e*x^4)^(n - 1), x)}}\\
{\leu0{ - ((U - V*v)*d + (V*c + 3*W*b)*w + 4*((U + n*V*d)*e + W*w*c)*x}}\\
{\leu1{ + 2*(V*e + 2*W*d)*w*x^2 + 4*W*w*e*x^3)}}\\
{\leu1{*(a + b*x + c*x^2 + d*x^3 + e*x^4)^n = 0}},
\bds\bfl
where {\ttt{U = 2*(2*n + 1)*((4*n + 3)*A*e - (2*n + 1)*C*c)*e}}, 
{\ttt{V = (4*n + 3)*B*e -{\fb} (3*n + 2)*C*d}}, 
{\ttt{W = (2*n + 1)*C*e}}, {\ttt{v = (3*n + 1)*d}}, 
{\ttt{w = 2*(4*n + 1)*e}}.
\efl

\subsection{Singly degenerate integrands}
\label{3B}

Here the quartic polynomial is required to possess one double root whence 
its discriminant must vanish: {\ttt{sa*re - 2*sb*rd +{\fb} sc*rc1 = 0}}. 
To exclude doubly degenerate quartics, the recurrences should be applied 
only if additionally {\ttt{sa \tneq 0}}.{\pagebreak[3]}
\ads
{\len0{u*(u + 1)*(sa*c - 3*sb*d + 6*sc*e)*sa*INT((A + B*x)}}\\
{\leu1{*(a + b*x + c*x^2 + d*x^3 + e*x^4)^n, x)}}\\
{\leu0{ - 2*INT((2*(4*n + 5)*(U*rb - V)*e - (3*n + 4)*(U*ra - W)*d}}\\
{\leu1{ - 2*(2*n + 3)*(U*ra - W)*e*x)}}\\
{\leu1{*(a + b*x + c*x^2 + d*x^3 + e*x^4)^(n + 1), x)}}\\
{\leu0{ - (u + 1)*(A*sa - B*sb)*(sb*c - 2*sc*d + 3*sd*e}}\\
{\leu1{ + (sa*c - sb*d + sc*e)*x + (sa*d - sb*e)*x^2 + sa*e*x^3)}}\\
{\leu1{*(a + b*x + c*x^2 + d*x^3 + e*x^4)^n}}\\
{\leu0{ + ((U*rb - V)*d - (U*ra - W)*c + 2*(2*(U*rb - V)*e - (U*ra - W)*d)*x}}\\
{\leu1{ - 2*(U*ra - W)*e*x^2)*(a + b*x + c*x^2 + d*x^3 + e*x^4)^(n + 1) = 0}},
\bds\bfl
where {\ttt{U = (3*n + 2)*A*sa*d - (2*n + 1)*(2*A*sb*e + B*sa*c) +%
{\fb} B*((3*n + 1)*sb*d - (4*n + 1)*sc*e)}}, 
{\ttt{V = (4*n + 3)*(2*A*rc2 - B*rd)*sa*e}}, 
{\ttt{W = (4*n + 3)*(2*A*rb - B*rc1)*sa*e}}, {\ttt{u = 2*n + 1}}.
\efl\ads
{\len0{8*(2*n + 1)*(4*n + 1)*w*e^2*INT((A + B*x)}}\\
{\leu1{*(a + b*x + c*x^2 + d*x^3 + e*x^4)^n, x)}}\\
{\leu0{ + n*INT(((U*ra + V*rb)*((2*n - 1)*sa*c - (3*n - 2)*sb*d}}\\
{\leu1{ + (4*n - 3)*sc*e) - (U*rc1 + V*rd)*w + ((U*ra + V*rb)}}\\
{\leu1{*((3*n - 1)*sa*d - 2*(2*n - 1)*sb*e) - 2*(U*rb + V*rc2)*w)*x)}}\\
{\leu1{*(a + b*x + c*x^2 + d*x^3 + e*x^4)^(n - 1), x)}}\\
{\leu0{ - n*(U*ra + V*rb)*(sb*c - 2*sc*d + 3*sd*e}}\\
{\leu1{ + (sa*c - sb*d + sc*e)*x + (sa*d - sb*e)*x^2 + sa*e*x^3)}}\\
{\leu1{*(a + b*x + c*x^2 + d*x^3 + e*x^4)^(n - 1)}}\\
{\leu0{ - w*(U*d + V*c + 2*(2*U*e + V*d)*x + 2*V*e*x^2)}}\\
{\leu1{*(a + b*x + c*x^2 + d*x^3 + e*x^4)^n = 0}},
\bds\bfl
where {\ttt{U = 2*(2*n + 1)*A*e - (3*n + 1)*B*d}}, 
{\ttt{V = 2*(4*n + 1)*B*e}}, {\ttt{w = (4*n - 1)*sa*e}}.
\efl

\subsection{Doubly degenerate integrands}
\label{3C}

Now the quartic polynomial is required to possess a triple root whence 
the following must vanish: {\ttt{6*ra*a - 3*rb*b + rc1*c ={\fb} %
rc1 - rc2 = 0}}. To exclude triply degenerate quartics, the recurrence 
should be applied only if additionally {\ttt{ra \tneq 0}}. Note that the 
relation can be used backwards for its own inverse.
\ads
{\len0{(n + 1)*(3*n + 1)*(3*n + 2)*ra^2}}\\
{\leu1{*INT((a + b*x + c*x^2 + d*x^3 + e*x^4)^n, x)}}\\
{\leu0{ + 4*(4*n + 5)*w*e*INT((a + b*x + c*x^2 + d*x^3 + e*x^4)^(n + 1), x)}}\\
{\leu0{ - (n + 1)*(3*(5*n + 3)*ra*b*e - v*c + ((22*n + 13)*ra*c*e - 3*v*d)*x}}\\
{\leu1{ + 3*((7*n + 4)*ra*d - 2*v)*e*x^2 + 6*(2*n + 1)*ra*e^2*x^3)}}\\
{\leu1{*(a + b*x + c*x^2 + d*x^3 + e*x^4)^n}}\\
{\leu0{ - w*(d + 4*e*x)*(a + b*x + c*x^2 + d*x^3 + e*x^4)^(n + 1) = 0}},
\bds\bfl
where {\ttt{v = (3*n + 2)*ra*d - (4*n + 3)*rb*e}}, 
{\ttt{w = 6*(2*n + 1)*(4*n + 3)*e^2}}.
\efl

Now the quartic polynomial is required to possess two double roots whence 
the following must vanish: {\ttt{ra*d - 4*rb*e = rb*d - 4*rc2*e = 0}}. To 
exclude triply degenerate quartics, the recurrence should be applied only 
if additionally {\ttt{ra \tneq 0}}. Note that the relation can be used 
backwards for its own inverse.
\ads
{\len0{2*(n + 1)*(2*n + 1)*(rc2 - rc1)}}\\
{\leu1{*INT((a + b*x + c*x^2 + d*x^3 + e*x^4)^n, x)}}\\
{\leu0{ + 4*w*(w + 2)*e*INT((a + b*x + c*x^2 + d*x^3 + e*x^4)^(n + 1), x)}}\\
{\leu0{ - (n + 1)*(rd + (rc1 + 2*rc2)*x + 3*rb*x^2 + ra*x^3)}}\\
{\leu1{*(a + b*x + c*x^2 + d*x^3 + e*x^4)^n}}\\
{\leu0{ - w*(d + 4*e*x)*(a + b*x + c*x^2 + d*x^3 + e*x^4)^(n + 1) = 0}}, 
\bds\bfl
where {\ttt{w = 4*n + 3}}.
\efl

\numberwithin{equation}{section}
\tttdef\formiv{(a + b*x)^m *(c + d*x)^n}
\section{Integrands involving \formiv}

To exclude integrands with confluent roots, the following recurrences 
should be applied only if the resultant of the linear polynomials does 
not vanish: {\ttt{a*d - b*c \tneq 0}}. Note that each relation can be 
used backwards for its own inverse.
\ads
{\len0{(m + 1)*(a*d - b*c)*INT((a + b*x)^m *(c + d*x)^n, x)}}\\
{\leu0{ - d*(m + n + 2)*INT((a + b*x)^(m + 1)*(c + d*x)^n, x)}}\\
{\leu0{ + (a + b*x)^(m + 1)*(c + d*x)^(n + 1) = 0}}.
\ads
{\len0{(m + 1)*b*INT((a + b*x)^m *(c + d*x)^n, x)}}\\
{\leu0{ + n*d*INT((a + b*x)^(m + 1)*(c + d*x)^(n - 1), x)}}\\
{\leu0{ - (a + b*x)^(m + 1)*(c + d*x)^n = 0}}.

\tttdef\formv{exp(a + b*x)*(c + d*x)^n}
\section{Integrands involving \formv}

Note that each recurrence can be used backwards for its own (or its 
partner's) inverse.
\ads
{\len0{(n + 1)*d*INT(exp(a + b*x)*(c + d*x)^n, x)}}\\
{\leu0{ + b*INT(exp(a + b*x)*(c + d*x)^(n + 1), x)}}\\
{\leu0{ - exp(a + b*x)*(c + d*x)^(n + 1) = 0}}.
\ads
{\len0{(n + 1)*d*INT(cos(a + b*x)*(c + d*x)^n, x)}}\\
{\leu0{ - b*INT(sin(a + b*x)*(c + d*x)^(n + 1), x)}}\\
{\leu0{ - cos(a + b*x)*(c + d*x)^(n + 1) = 0}}.
\ads
{\len0{(n + 1)*d*INT(sin(a + b*x)*(c + d*x)^n, x)}}\\
{\leu0{ + b*INT(cos(a + b*x)*(c + d*x)^(n + 1), x)}}\\
{\leu0{ - sin(a + b*x)*(c + d*x)^(n + 1) = 0}}.

\tttdef\formvi{(a + b*x)^m *(c + d*x + e*x^2)^n}
\section{Integrands involving \formvi}
\numberwithin{equation}{subsection}

Abbreviations used with these integrands:
\ads
{\reu0{ra = 2*a*e - b*d}}, {\ttt{rb = a*d - 2*b*c}}.

\subsection{Nondegenerate integrands}
\label{6A}

To exclude integrands with confluent roots, the following recurrences 
should be applied only if the overall discriminant does not vanish: 
{\ttt{(ra*a - rb*b)*(4*c*e - d^2) \tneq 0}}.
\ads
{\len0{(m + 1)*(ra*a - rb*b)*INT((A + B*x)}}\\
{\leu1{*(a + b*x)^m *(c + d*x + e*x^2)^n, x)}}\\
{\leu0{ - INT(((m + 1)*(A*ra - B*rb) - v*V*d - 2*v*V*e*x)}}\\
{\leu1{*(a + b*x)^(m + 1)*(c + d*x + e*x^2)^n, x)}}\\
{\leu0{ - 2*V*(a + b*x)^(m + 1)*(c + d*x + e*x^2)^(n + 1) = 0}},
\bds\bfl
where {\ttt{V = A*b - B*a}}, {\ttt{v = m + 2*n + 3}}.
\efl\ads
{\len0{2*u*e*INT((A + B*x)*(a + b*x)^m *(c + d*x + e*x^2)^n, x)}}\\
{\leu0{ - INT((u*U*a + m*B*rb + (u*U*b + m*B*ra)*x)}}\\
{\leu1{*(a + b*x)^(m - 1)*(c + d*x + e*x^2)^n, x)}}\\
{\leu0{ - 2*B*(a + b*x)^m *(c + d*x + e*x^2)^(n + 1) = 0}},
\bds\bfl
where {\ttt{U = 2*A*e - B*d}}, {\ttt{u = m + 2*n + 2}}.
\efl\ads
{\len0{(n + 1)* 1/2 *(ra*a - rb*b)*(4*c*e - d^2)*INT((A + B*x)}}\\
{\leu1{*(a + b*x)^m *(c + d*x + e*x^2)^n, x)}}\\
{\leu0{ + INT((U*u - V*((n + 1)*ra + a*e) - V*(u + b)*e*x)}}\\
{\leu1{*(a + b*x)^m *(c + d*x + e*x^2)^(n + 1), x)}}\\
{\leu0{ - (U - V*e*x)*(a + b*x)^(m + 1)*(c + d*x + e*x^2)^(n + 1) = 0}},
\bds\bfl
where {\ttt{U = A*rb*e - (A*d - B*c)*ra}}, {\ttt{V = A*ra - B*rb}}, 
{\ttt{u = (m + 2*n + 3)*b}}.
\efl\ads
{\len0{2*v*(v + b)*e*INT((A + B*x)*(a + b*x)^m *(c + d*x + e*x^2)^n, x)}}\\
{\leu0{ + n*INT((U*rb - V*v*a + (U*ra - V*v*b)*x)}}\\
{\leu1{*(a + b*x)^m *(c + d*x + e*x^2)^(n - 1), x)}}\\
{\leu0{ - (U + B*v*d + 2*B*v*e*x)*(a + b*x)^(m + 1)*(c + d*x + e*x^2)^n = 0}},
\bds\bfl
where {\ttt{U = 2*(m + 2*n + 2)*(A*b - B*a)*e + (m + 1)*B*ra}}, 
{\ttt{V = B*(4*c*e - d^2)}}, {\ttt{v = (m + 2*n + 1)*b}}.
\efl\ads
{\len0{(n + 1)*(4*c*e - d^2)*INT((A + B*x)}}\\
{\leu1{*(a + b*x)^m *(c + d*x + e*x^2)^n, x)}}\\
{\leu0{ - INT((m*U*b + v*V*a + (v + m)*V*b*x)}}\\
{\leu1{*(a + b*x)^(m - 1)*(c + d*x + e*x^2)^(n + 1), x)}}\\
{\leu0{ + (U + V*x)*(a + b*x)^m *(c + d*x + e*x^2)^(n + 1) = 0}},
\bds\bfl
where {\ttt{U = A*d - 2*B*c}}, {\ttt{V = 2*A*e - B*d}}, 
{\ttt{v = 2*n + 3}}.
\efl\ads
{\len0{(m + 1)*u*b^2*INT((A + B*x)*(a + b*x)^m *(c + d*x + e*x^2)^n, x)}}\\
{\leu0{ + n*INT((u*U*d + V*rb + (2*u*U*e + V*ra)*x)}}\\
{\leu1{*(a + b*x)^(m + 1)*(c + d*x + e*x^2)^(n - 1), x)}}\\
{\leu0{ - (u*U + V*a + V*b*x)*(a + b*x)^(m + 1)*(c + d*x + e*x^2)^n = 0}},
\bds\bfl
where {\ttt{U = A*b - B*a}}, {\ttt{V = (m + 1)*B}}, 
{\ttt{u = m + 2*n + 2}}.
\efl

\subsection{Singly degenerate integrands}
\label{6B}

Now the root of the linear factor is required to coincide with one root 
of the quadratic whence the corresponding resultant must vanish: 
{\ttt{ra*a - rb*b = 0}}. To exclude additional degeneracies, the 
recurrences should be applied only if also {\ttt{ra \tneq 0}}. Note that 
each relation can be used backwards for its own inverse.
\ads
{\len0{(m + n + 1)*ra*INT((a + b*x)^m *(c + d*x + e*x^2)^n, x)}}\\
{\leu0{ - (m + 2*n + 2)*e*INT((a + b*x)^(m + 1)*(c + d*x + e*x^2)^n, x)}}\\
{\leu0{ + b*(a + b*x)^m *(c + d*x + e*x^2)^(n + 1) = 0}}.
\ads
{\len0{(m + n + 1)*v*(4*c*e - d^2)*INT((a + b*x)^m *(c + d*x + e*x^2)^n, x)}}\\
{\leu0{ - w*(w + 1)*b*e*INT((a + b*x)^m *(c + d*x + e*x^2)^(n + 1), x)}}\\
{\leu0{ + (m*a*e + v*d + w*b*e*x)*(a + b*x)^m *(c + d*x + e*x^2)^(n + 1) = 0}},
\bds\bfl
where {\ttt{v = (n + 1)*b}}, {\ttt{w = m + 2*n + 2}}.
\efl\ads
{\len0{(n + 1)*ra*INT((a + b*x)^m *(c + d*x + e*x^2)^n, x)}}\\
{\leu0{ + (m + 2*n + 2)*b^2*INT((a + b*x)^(m - 1)*(c + d*x + e*x^2)^(n + 1), x)}}\\
{\leu0{ - b*(a + b*x)^m *(c + d*x + e*x^2)^(n + 1) = 0}}.
\bds

Now the quadratic polynomial is required to possess a double root whence 
its discriminant must vanish: {\ttt{4*c*e - d^2 = 0}}. To exclude 
additional degeneracies, the recurrences should be applied only if also 
{\ttt{ra \tneq 0}}. Note that each relation can be used backwards for its 
own inverse.
\ads
{\len0{(m + 1)*ra*INT((a + b*x)^m *(c + d*x + e*x^2)^n, x)}}\\
{\leu0{ - 2*(m + 2*n + 2)*e*INT((a + b*x)^(m + 1)*(c + d*x + e*x^2)^n, x)}}\\
{\leu0{ + (d + 2*e*x)*(a + b*x)^(m + 1)*(c + d*x + e*x^2)^n = 0}}.
\ads
{\len0{(n + 1)*(2*n + 1)*(ra*a - rb*b)\\}}
{\leu1{*INT((a + b*x)^m *(c + d*x + e*x^2)^n, x)}}\\
{\leu0{ - w*(w + b)*INT((a + b*x)^m *(c + d*x + e*x^2)^(n + 1), x)}}\\
{\leu0{ - (n + 1)*(rb + ra*x)*(a + b*x)^(m + 1)*(c + d*x + e*x^2)^n}}\\
{\leu0{ + w*(a + b*x)^(m + 1)*(c + d*x + e*x^2)^(n + 1) = 0}},
\bds\bfl
where {\ttt{w = (m + 2*n + 2)*b}}.
\efl\ads
{\len0{(n + 1)*(2*n + 1)*ra*INT((a + b*x)^m *(c + d*x + e*x^2)^n, x)}}\\
{\leu0{ - m*(m + 2*n + 2)*b^2}}\\
{\leu1{*INT((a + b*x)^(m - 1)*(c + d*x + e*x^2)^(n + 1), x)}}\\
{\leu0{ - (n + 1)*(rb + ra*x)*(a + b*x)^m *(c + d*x + e*x^2)^n}}\\
{\leu0{ + m*b*(a + b*x)^m *(c + d*x + e*x^2)^(n + 1) = 0}}.

\tttdef\formvii{exp(a + b*x)*(c + d*x + e*x^2)^n}
\section{Integrands involving \formvii}

\subsection{Nondegenerate integrands}
\label{7A}

To exclude integrands with confluent roots, the following recurrences 
should be applied only if the quadratic discriminant does not vanish: 
{\ttt{4*c*e - d^2 \tneq 0}}.
\ads
{\len0{(n + 1)*(4*c*e - d^2)*INT((A + B*x)}}\\
{\leu1{*exp(a + b*x)*(c + d*x + e*x^2)^n, x)}}\\
{\leu0{ - INT((U*b + (2*n + 3)*V + V*b*x)}}\\
{\leu1{*exp(a + b*x)*(c + d*x + e*x^2)^(n + 1), x)}}\\
{\leu0{ + (U + V*x)*exp(a + b*x)*(c + d*x + e*x^2)^(n + 1) = 0}},
\bds\bfl
where {\ttt{U = A*d - 2*B*c}}, {\ttt{V = 2*A*e - B*d}}.
\efl\ads
{\len0{b^2*INT((A + B*x)*exp(a + b*x)*(c + d*x + e*x^2)^n, x)}}\\
{\leu0{ + n*INT((U*d - 2*B*b*c + (2*U*e - B*b*d)*x)}}\\
{\leu1{*exp(a + b*x)*(c + d*x + e*x^2)^(n - 1), x)}}\\
{\leu0{ - (U + B*b*x)*exp(a + b*x)*(c + d*x + e*x^2)^n = 0}},
\bds\bfl
where {\ttt{U = A*b - (2*n + 1)*B}}.
\efl

\subsection{Degenerate integrands}
\label{7B}

Now the quadratic polynomial is required to possess a double root whence 
its discriminant must vanish: {\ttt{4*c*e - d^2 = 0}}. Note that each 
recurrence can be used backwards for its own inverse.
\ads
{\len0{2*(n + 1)*(2*n + 1)*e*INT(exp(a + b*x)*(c + d*x + e*x^2)^n, x)}}\\
{\leu0{ - b^2*INT(exp(a + b*x)*(c + d*x + e*x^2)^(n + 1), x)}}\\
{\leu0{ - (n + 1)*(d + 2*e*x)*exp(a + b*x)*(c + d*x + e*x^2)^n}}\\
{\leu0{ + b*exp(a + b*x)*(c + d*x + e*x^2)^(n + 1) = 0}}.
\ads
{\len0{2*(n + 1)*(2*n + 1)*e*INT(cos(a + b*x)*(c + d*x + e*x^2)^n, x)}}\\
{\leu0{ + b^2*INT(cos(a + b*x)*(c + d*x + e*x^2)^(n + 1), x)}}\\
{\leu0{ - (n + 1)*(d + 2*e*x)*cos(a + b*x)*(c + d*x + e*x^2)^n}}\\
{\leu0{ - b*sin(a + b*x)*(c + d*x + e*x^2)^(n + 1) = 0}}.
\ads
{\len0{2*(n + 1)*(2*n + 1)*e*INT(sin(a + b*x)*(c + d*x + e*x^2)^n, x)}}\\
{\leu0{ + b^2*INT(sin(a + b*x)*(c + d*x + e*x^2)^(n + 1), x)}}\\
{\leu0{ - (n + 1)*(d + 2*e*x)*sin(a + b*x)*(c + d*x + e*x^2)^n}}\\
{\leu0{ + b*cos(a + b*x)*(c + d*x + e*x^2)^(n + 1) = 0}}.

\tttdef\formviii{(a + b*x + c*x^2)^m *(d + e*x + f*x^2)^n}
\section{Integrands involving \formviii}

Abbreviations used with these integrands:
\ads
{\reu0{ra = b*f - c*e}}, {\ttt{rb = a*f - c*d}}, {\ttt{rc = a*e - b*d}}.\\[1ex]
{\reu0{sa = ra*b - 2*rb*c}}, {\ttt{sb = ra*a - rc*c}}, 
{\ttt{sc = 2*rb*a - rc*b}},\\
{\reu0{sd = 2*rb*f - ra*e}}, {\ttt{se = rc*f - ra*d}}.

\subsection{Nondegenerate integrands}
\label{8A}

To exclude integrands with confluent roots, the following recurrences 
should be applied only if the overall discriminant does not vanish: 
{\ttt{(ra*rc - rb^2)*(4*a*c - b^2)*(4*d*f - e^2) \tneq 0}}.
\ads
{\len0{u*(4*a*c - b^2)*INT((A + B*x + C*x^2)}}\\
{\leu1{*(a + b*x + c*x^2)^m *(d + e*x + f*x^2)^n, x)}}\\
{\leu0{ - INT(((2*A*c - B*b + 2*C*a)*u - v*V*e + W*((m + 1)*rb - c*d)}}\\
{\leu1{ - (2*v*V*f - W*((m + 1)*ra - (v + 1)*c*e))*x - (2*v + 1)*W*c*f*x^2)}}\\
{\leu1{*(a + b*x + c*x^2)^(m + 1)*(d + e*x + f*x^2)^n, x)}}\\
{\leu0{ - (V + W*c*x)*(a + b*x + c*x^2)^(m + 1)*(d + e*x + f*x^2)^(n + 1) = 0}},
\bds\bfl
where {\ttt{V = (A*b - B*a)*sa - (A*c - C*a)*sb}}, 
{\ttt{W = A*sa - B*sb + C*sc}}, {\ttt{u = (m + 1)*(ra*rc - rb^2)}}, 
{\ttt{v = m + n + 2}}.
\efl\ads
{\len0{2*w*(2*w + f)*c*INT((A + B*x + C*x^2)}}\\
{\leu1{*(a + b*x + c*x^2)^m *(d + e*x + f*x^2)^n, x)}}\\
{\leu0{ - INT((U*a + m*V*rc + W*w*d + (U*b + 2*m*V*rb + W*w*e)*x}}\\
{\leu1{ + (U*c + m*V*ra + W*w*f)*x^2)}}\\
{\leu1{*(a + b*x + c*x^2)^(m - 1)*(d + e*x + f*x^2)^n, x)}}\\
{\leu0{ - (V + C*w*b + 2*C*w*c*x)}}\\
{\leu1{*(a + b*x + c*x^2)^m *(d + e*x + f*x^2)^(n + 1) = 0}},
\bds\bfl
where {\ttt{U = (m + n + 1)*((2*m + 2*n + 3)*(2*A*f - B*e + 2*C*d)*f -%
{\fb} (2*m + n + 2)*C*(4*d*f - e^2))*c}}, 
{\ttt{V = (2*m + 2*n + 3)*(B*f - C*e)*c -{\fb} (n + 1)*C*ra}}, 
{\ttt{W = m*C*(4*a*c - b^2)}}, {\ttt{w = (m + n + 1)*f}}.
\efl\ads
{\len0{u*c*INT((A + B*x + C*x^2)*(a + b*x + c*x^2)^m *(d + e*x + f*x^2)^n, x)}}\\
{\leu0{ - INT((n*V*e + (C*u + w*W)*d + (2*n*V*f + (C*u + (w + n)*W)*e)*x}}\\
{\leu1{ + (C*u + (w + 2*n)*W)*f*x^2)}}\\
{\leu1{*(a + b*x + c*x^2)^(m + 1)*(d + e*x + f*x^2)^(n - 1), x)}}\\
{\leu0{ + (V + W*x)*(a + b*x + c*x^2)^(m + 1)*(d + e*x + f*x^2)^n = 0}},
\bds\bfl
where {\ttt{V = A*b*c - 2*B*a*c + C*a*b}}, 
{\ttt{W = 2*(A*c - C*a)*c - (B*c - C*b)*b}}, 
{\ttt{u = (m + 1)*(4*a*c - b^2)}}, {\ttt{w = 2*m + 3}}.
\efl

\subsection{Integrands with a single cross-degeneracy}
\label{8B}

Here the two quadratics are required to have one root in common whence 
their resultant must vanish: {\ttt{ra*rc -{\fb} rb^2 = 0}}. To exclude 
additional degeneracies, the recurrences should be applied only if also 
{\ttt{(2*a*f - b*e + 2*c*d)*ra \tneq 0}}.
\ads
{\len0{(m + n + 1)*v*sd*(4*a*c - b^2)*INT((A + B*x)}}\\
{\leu1{*(a + b*x + c*x^2)^m *(d + e*x + f*x^2)^n, x)}}\\
{\leu0{ + INT((V*(u*c*e - v) + W*(u*b*f + n*ra) + 2*u*(V + W)*c*f*x)}}\\
{\leu1{*(a + b*x + c*x^2)^(m + 1)*(d + e*x + f*x^2)^n, x)}}\\
{\leu0{ + ra*(n*A*se*c + ((m + 1)*A*b - (2*m + n + 2)*B*a)*sd + V*x)}}\\
{\leu1{*(a + b*x + c*x^2)^(m + 1)*(d + e*x + f*x^2)^n}}\\
{\leu0{ - (V + W)*f*(a + b*x + c*x^2)^(m + 2)*(d + e*x + f*x^2)^n = 0}},
\bds\bfl
where {\ttt{V = ((2*m + n + 2)*A*c - (m + n + 1)*B*b)*sd + n*B*se*c}}, 
{\ttt{W = (m + 1)*(A*sa - B*sb)*f}}, {\ttt{u = m + n + 2}}, 
{\ttt{v = (m + 1)*ra}}.
\efl\ads
{\len0{w*(w + f)*ra*c*INT((A + B*x)}}\\
{\leu1{*(a + b*x + c*x^2)^m *(d + e*x + f*x^2)^n, x)}}\\
{\leu0{ - INT((U*((m + n)*sd*b - n*se*c) + m*V*sb}}\\
{\leu1{ + ((2*m + n)*U*sd*c + m*V*sa)*x)}}\\
{\leu1{*(a + b*x + c*x^2)^(m - 1)*(d + e*x + f*x^2)^n, x)}}\\
{\leu0{ + (U*sd - V*rb - V*ra*x)*(a + b*x + c*x^2)^m *(d + e*x + f*x^2)^n}}\\
{\leu0{ - B*w*ra*f*(a + b*x + c*x^2)^(m + 1)*(d + e*x + f*x^2)^n = 0}},
\bds\bfl
where {\ttt{U = (m + n + 1)*(2*A*f - B*e)*c + m*B*ra}}, 
{\ttt{V = ((m + n + 1)*(2*A*c - B*b)*f - n*B*ra)*f}}, 
{\ttt{w = (2*m + 2*n + 1)*f}}.
\efl\ads
{\len0{2*(m + 1)*sa*c*INT((A + B*x)}}\\
{\leu1{*(a + b*x + c*x^2)^m *(d + e*x + f*x^2)^n, x)}}\\
{\leu0{ + INT(((m + 1)*U*sd - V*e - 2*V*f*x)}}\\
{\leu1{*(a + b*x + c*x^2)^(m + 1)*(d + e*x + f*x^2)^(n - 1), x)}}\\
{\leu0{ + (U*ra - B*sa)*(a + b*x + c*x^2)^(m + 1)*(d + e*x + f*x^2)^n = 0}},
\bds\bfl
where {\ttt{U = 2*A*c - B*b}}, 
{\ttt{V = (m + n + 1)*(2*A*c - B*b)*ra - n*B*sa}}.
\efl

\subsection{Integrands with a single self-degeneracy}
\label{8C}

Here the first quadratic polynomial is required to possess a double root 
whence its discriminant must vanish: {\ttt{4*a*c - b^2 = 0}}. To exclude 
additional degeneracies, the recurrences should be applied only if also 
{\ttt{(2*a*f - b*e + 2*c*d)*(4*d*f - e^2) \tneq 0}}.
\ads
{\len0{(m + 1)*u*v^2*INT((A + B*x)}}\\
{\leu1{*(a + b*x + c*x^2)^m *(d + e*x + f*x^2)^n, x)}}\\
{\leu0{ + INT(((u + 1)*((n + 1)*A*v - 2*(u + n + 1)*(A*rb - B*rc))*f}}\\
{\leu1{ + (u + n + 2)*W*e + 2*(m + n + 2)*W*f*x)}}\\
{\leu1{*(a + b*x + c*x^2)^(m + 1)*(d + e*x + f*x^2)^n, x)}}\\
{\leu0{ - (m + 1)*v*(A*b - 2*B*a + (2*A*c - B*b)*x)}}\\
{\leu1{*(a + b*x + c*x^2)^m *(d + e*x + f*x^2)^(n + 1)}}\\
{\leu0{ - W*(a + b*x + c*x^2)^(m + 1)*(d + e*x + f*x^2)^(n + 1) = 0}},
\bds\bfl
where {\ttt{W = 2*(2*m + n + 2)*(A*ra - B*rb) - (n + 1)*B*%
(2*a*f - b*e + 2*c*d)}}, {\ttt{u = 2*m + 1}}, 
{\ttt{v = 2*a*f - b*e + 2*c*d}}.
\efl\ads
{\len0{2*w*(w + f)*INT((A + B*x)*(a + b*x + c*x^2)^m *(d + e*x + f*x^2)^n, x)}}\\
{\leu0{ - INT((U*(u*rb + v) + u*V*rc + (u*U*ra + V*(u*rb - v))*x)}}\\
{\leu1{*(a + b*x + c*x^2)^(m - 1)*(d + e*x + f*x^2)^n, x)}}\\
{\leu0{ - (U*b + 2*V*a + (2*U*c + V*b)*x)}}\\
{\leu1{*(a + b*x + c*x^2)^(m - 1)*(d + e*x + f*x^2)^(n + 1)}}\\
{\leu0{ - 2*B*w*(a + b*x + c*x^2)^m *(d + e*x + f*x^2)^(n + 1) = 0}},
\bds\bfl
where {\ttt{U = (m + n + 1)*(2*A*f - B*e) - m*B*e}}, {\ttt{V = 2*m*B*f}}, 
{\ttt{u = 2*(2*m + n)}}, {\ttt{v = (n + 1)*(2*a*f - b*e + 2*c*d)}}, 
{\ttt{w = (2*m + 2*n + 1)*f}}.
\efl\ads
{\len0{(n + 1)*u*v*INT((A + B*x)*(a + b*x + c*x^2)^m *(d + e*x + f*x^2)^n, x)}}\\
{\leu0{ - INT((m*(2*A*c - B*b)*u + (2*n + 3)*V*v - w*W*e - 2*w*W*f*x)}}\\
{\leu1{*(a + b*x + c*x^2)^m *(d + e*x + f*x^2)^(n + 1), x)}}\\
{\leu0{ + v*(A*e - 2*B*d + V*x)*(a + b*x + c*x^2)^m *(d + e*x + f*x^2)^(n + 1)}}\\
{\leu0{ - W*(a + b*x + c*x^2)^m *(d + e*x + f*x^2)^(n + 2) = 0}},
\bds\bfl
where {\ttt{V = 2*A*f - B*e}}, 
{\ttt{W = 2*(A*e - 2*B*d)*c - (2*A*f - B*e)*b}}, {\ttt{u = 4*d*f - e^2}}, 
{\ttt{v = 2*a*f - b*e + 2*c*d}}, {\ttt{w = m + n + 2}}.
\efl\ads
{\len0{2*w*(w + 1)*c*f*INT((A + B*x)}}\\
{\leu1{*(a + b*x + c*x^2)^m *(d + e*x + f*x^2)^n, x)}}\\
{\leu0{ + n*INT((U*e - 2*V*d + (2*U*f - V*e)*x)}}\\
{\leu1{*(a + b*x + c*x^2)^m *(d + e*x + f*x^2)^(n - 1), x)}}\\
{\leu0{ - (U + V*x)*(a + b*x + c*x^2)^m *(d + e*x + f*x^2)^n}}\\
{\leu0{ - 2*B*w*c*(a + b*x + c*x^2)^m *(d + e*x + f*x^2)^(n + 1) = 0}},
\bds\bfl
where {\ttt{U = (m + n + 1)*(2*A*f - B*e)*b +%
{\fb} B*(2*m*rb - (m + 2*n + 1)*{\fb}(2*a*f - b*e + 2*c*d))}}, 
{\ttt{V = 2*((m + n + 1)*(2*A*f - B*e)*c + m*B*ra)}}, 
{\ttt{w = 2*m + 2*n + 1}}.
\efl\ads
{\len0{2*(m + 1)*(2*m + 1)*(2*a*f - b*e + 2*c*d)*c*INT((A + B*x)}}\\
{\leu1{*(a + b*x + c*x^2)^m *(d + e*x + f*x^2)^n, x)}}\\
{\leu0{ + n*INT((W*(4*d*f - e^2) + V*e + 2*V*f*x)}}\\
{\leu1{*(a + b*x + c*x^2)^(m + 1)*(d + e*x + f*x^2)^(n - 1), x)}}\\
{\leu0{ - (V - W*e - 2*W*f*x)*(a + b*x + c*x^2)^(m + 1)*(d + e*x + f*x^2)^n}}\\
{\leu0{ - W*(b + 2*c*x)*(a + b*x + c*x^2)^m *(d + e*x + f*x^2)^(n + 1) = 0}},
\bds\bfl
where {\ttt{V = 2*(m + n + 1)*(A*ra - B*rb) + (m - n)*B*%
(2*a*f - b*e + 2*c*d)}}, {\ttt{W = (m + 1)*(2*A*c - B*b)}}.
\efl\ads
{\len0{2*(n + 1)*u*f*INT((A + B*x)}}\\
{\leu1{*(a + b*x + c*x^2)^m *(d + e*x + f*x^2)^n, x)}}\\
{\leu0{ + INT((m*B*u*b - V*(v*rb - (m - n - 2)*(2*a*f - b*e + 2*c*d))}}\\
{\leu1{ + (2*m*B*u*c - v*V*ra)*x)}}\\
{\leu1{*(a + b*x + c*x^2)^(m - 1)*(d + e*x + f*x^2)^(n + 1), x)}}\\
{\leu0{ + 2*f*(A*e - 2*B*d + V*x)}}\\
{\leu1{*(a + b*x + c*x^2)^m *(d + e*x + f*x^2)^(n + 1)}}\\
{\leu0{ - V*(b + 2*c*x)*(a + b*x + c*x^2)^(m - 1)*(d + e*x + f*x^2)^(n + 2) = 0}},
\bds\bfl
where {\ttt{V = 2*A*f - B*e}}, {\ttt{u = 4*d*f - e^2}}, 
{\ttt{v = 2*(m + n + 1)}}.
\efl

\subsection{Doubly degenerate integrands}
\label{8D}

Now the two quadratics are required to have two roots in common whence 
the following must vanish: {\ttt{ra = rb = 0}}. To exclude additional 
degeneracies, the recurrences should be applied only if also 
{{\ttt2*a*f - b*e + 2*c*d \tneq 0}}. Note that each relation can be used 
backwards for its own inverse.
\ads
{\len0{u*(4*a*c - b^2)*INT((a + b*x + c*x^2)^m *(d + e*x + f*x^2)^n, x)}}\\
{\leu0{ - 2*(2*u + 1)*c*INT((a + b*x + c*x^2)^(m + 1)*(d + e*x + f*x^2)^n, x)}}\\
{\leu0{ + (b + 2*c*x)*(a + b*x + c*x^2)^(m + 1)*(d + e*x + f*x^2)^n = 0}},
\bds\bfl
where {\ttt{u = m + n + 1}}.
\efl\ads
{\len0{c*INT((a + b*x + c*x^2)^m *(d + e*x + f*x^2)^n, x)}}\\
{\leu0{ - f*INT((a + b*x + c*x^2)^(m + 1)*(d + e*x + f*x^2)^(n - 1), x) = 0}}.
\bds

Now the first quadratic is required to possess a double root which 
coincides with one root of the second quadratic whence the following must 
vanish: {\ttt{4*a*c - b^2 ={\fb} 2*a*f - b*e + 2*c*d = 0}}. To exclude 
additional degeneracies, the recurrences should be applied only if also 
{\ttt{4*d*f - e^2 \tneq 0}}. Note that each relation can be used 
backwards for its own inverse.
\ads
{\len0{(2*m + n + 1)*v*(4*d*f - e^2)*c}}\\
{\leu1{*INT((a + b*x + c*x^2)^m *(d + e*x + f*x^2)^n, x)}}\\
{\leu0{ + 2*w*(2*w + f)*INT((a + b*x + c*x^2)^(m + 1)*(d + e*x + f*x^2)^n, x)}}\\
{\leu0{ - (v*ra + w*b + 2*w*c*x)}}\\
{\leu1{*(a + b*x + c*x^2)^m *(d + e*x + f*x^2)^(n + 1) = 0}},
\bds\bfl
where {\ttt{v = 2*m + n + 2}}, {\ttt{w = (m + n + 1)*f}}.
\efl\ads
{\len0{(2*m + n + 1)*(n + 1)*(4*d*f - e^2)*c}}\\
{\leu1{*INT((a + b*x + c*x^2)^m *(d + e*x + f*x^2)^n, x)}}\\
{\leu0{ - 2*w*(2*w + 1)*c*f}}\\
{\leu1{*INT((a + b*x + c*x^2)^m *(d + e*x + f*x^2)^(n + 1), x)}}\\
{\leu0{ + (m*ra + w*c*e + 2*w*c*f*x)}}\\
{\leu1{*(a + b*x + c*x^2)^m *(d + e*x + f*x^2)^(n + 1) = 0}},
\bds\bfl
where {\ttt{w = m + n + 1}}.
\efl\ads
{\len0{2*(2*m + n + 1)*c*INT((a + b*x + c*x^2)^m *(d + e*x + f*x^2)^n, x)}}\\
{\leu0{ + 2*n*f*INT((a + b*x + c*x^2)^(m + 1)*(d + e*x + f*x^2)^(n - 1), x)}}\\
{\leu0{ - (b + 2*c*x)*(a + b*x + c*x^2)^m *(d + e*x + f*x^2)^n = 0}}.
\bds

Now both quadratics are required to possess a double root whence the 
following must vanish: {\ttt{4*a*c - b^2 = 4*d*f - e^2 = 0}}. To exclude 
additional degeneracies, the recurrences should be applied only if also 
{\ttt{ra \tneq 0}}. Note that each relation can be used backwards for its 
own inverse.
\ads
{\len0{2*(2*m + 1)*w*c*INT((a + b*x + c*x^2)^m *(d + e*x + f*x^2)^n, x)}}\\
{\leu0{ - 4*v*(2*v + 1)*c*f}}\\
{\leu1{*INT((a + b*x + c*x^2)^(m + 1)*(d + e*x + f*x^2)^n, x)}}\\
{\leu0{ - 2*(n*ra - v*b*f - 2*v*c*f*x)}}\\
{\leu1{*(a + b*x + c*x^2)^(m + 1)*(d + e*x + f*x^2)^n}}\\
{\leu0{ - w*(b + 2*c*x)*(a + b*x + c*x^2)^m *(d + e*x + f*x^2)^n = 0}},
\bds\bfl
where {\ttt{v = m + n + 1}}, {\ttt{w = (m + 1)*(2*a*f - b*e + 2*c*d)}}.
\efl\ads
{\len0{2*(m + 1)*(2*m + 1)*c*INT((a + b*x + c*x^2)^m *(d + e*x + f*x^2)^n, x)}}\\
{\leu0{ - 2*n*(2*n - 1)*f}}\\
{\leu1{*INT((a + b*x + c*x^2)^(m + 1)*(d + e*x + f*x^2)^(n - 1), x)}}\\
{\leu0{ - (m + 1)*(b + 2*c*x)*(a + b*x + c*x^2)^m *(d + e*x + f*x^2)^n}}\\
{\leu0{ + n*(e + 2*f*x)*(a + b*x + c*x^2)^(m + 1)*(d + e*x + f*x^2)^(n - 1) = 0}}.

\tttdef\formix{(a + b*x)^m *(c + d*x)^n *(e + f*x)^p}
\section{Integrands involving \formix}

Abbreviations used with these integrands:
\ads
{\reu0{ra = (c*f + d*e)*b - a*d*f}}, {\ttt{rb = b*c*e}}.

\subsection{Nondegenerate integrands}
\label{9A}

To exclude integrands with confluent roots, the following recurrences 
should be applied only if the overall discriminant does not vanish: 
{\ttt{(a*d - b*c)*(a*f - b*e)*{\fb}(c*f - d*e) \tneq 0}}.
\ads
{\len0\label{9A1}{(m + 1)*(a*d - b*c)*(a*f - b*e)*INT((A + B*x)}}\\
{\leu1{*(a + b*x)^m *(c + d*x)^n *(e + f*x)^p, x)}}\\
{\leu0{ + INT(((m + 1)*(A*ra - B*rb) + V*((n + 1)*d*e + (p + 1)*c*f)}}\\
{\leu1{ + (m + n + p + 3)*V*d*f*x)}}\\
{\leu1{*(a + b*x)^(m + 1)*(c + d*x)^n *(e + f*x)^p, x)}}\\
{\leu0{ - V*(a + b*x)^(m + 1)*(c + d*x)^(n + 1)*(e + f*x)^(p + 1) = 0}},
\bds\bfl
where {\ttt{V = A*b - B*a}}.
\efl\ads
{\len0{(m + n + p + 2)*d*f*INT((A + B*x)}}\\
{\leu1{*(a + b*x)^m *(c + d*x)^n *(e + f*x)^p, x)}}\\
{\leu0{ - INT((U*a - m*B*rb + (U*b - m*B*ra)*x)}}\\
{\leu1{*(a + b*x)^(m - 1)*(c + d*x)^n *(e + f*x)^p, x)}}\\
{\leu0{ - B*(a + b*x)^m *(c + d*x)^(n + 1)*(e + f*x)^(p + 1) = 0}},
\bds\bfl
where {\ttt{U = (m + n + p + 2)*A*d*f - B*((n + 1)*d*e + (p + 1)*c*f)}}.
\efl\ads
{\len0{u*b*INT((A + B*x)*(a + b*x)^m *(c + d*x)^n *(e + f*x)^p, x)}}\\
{\leu0{ - INT((B*u*c + V*(v*c*f + n*d*e) + (B*u + (v + n)*V*f)*d*x)}}\\
{\leu1{*(a + b*x)^(m + 1)*(c + d*x)^(n - 1)*(e + f*x)^p, x)}}\\
{\leu0{ + V*(a + b*x)^(m + 1)*(c + d*x)^n *(e + f*x)^(p + 1) = 0}},
\bds\bfl
where {\ttt{V = A*b - B*a}}, {\ttt{u = (m + 1)*(a*f - b*e)}}, 
{\ttt{v = m + p + 2}}.
\efl

\subsection{Singly degenerate integrands}
\label{9B}

Now the roots of the first two linear factors are required to coincide 
whence the corresponding resultant must vanish: {\ttt{a*d - b*c = 0}}. To 
exclude additional degeneracies, the recurrences should be applied only 
if also {\ttt{a*f - b*e \tneq 0}}. Note that each relation can be used 
backwards for its own inverse.
\ads
{\len0{(m + n + 1)*(a*f - b*e)*INT((a + b*x)^m *(c + d*x)^n *(e + f*x)^p, x)}}\\
{\leu0{ - (m + n + p + 2)*f*INT((a + b*x)^(m + 1)*(c + d*x)^n *(e + f*x)^p, x)}}\\
{\leu0{ + (a + b*x)^(m + 1)*(c + d*x)^n *(e + f*x)^(p + 1) = 0}}.
\ads
{\len0{(p + 1)*(a*f - b*e)*INT((a + b*x)^m *(c + d*x)^n *(e + f*x)^p, x)}}\\
{\leu0{ + (m + n + p + 2)*b*INT((a + b*x)^m *(c + d*x)^n *(e + f*x)^(p + 1), x)}}\\
{\leu0{ - (a + b*x)^(m + 1)*(c + d*x)^n *(e + f*x)^(p + 1) = 0}}.
\ads
{\len0{b*INT((a + b*x)^m *(c + d*x)^n *(e + f*x)^p, x)}}\\
{\leu0{ - d*INT((a + b*x)^(m + 1)*(c + d*x)^(n - 1)*(e + f*x)^p, x) = 0}}.
\ads
{\len0{(m + n + 1)*b*INT((a + b*x)^m *(c + d*x)^n *(e + f*x)^p, x)}}\\
{\leu0{ + p*f*INT((a + b*x)^(m + 1)*(c + d*x)^n *(e + f*x)^(p - 1), x)}}\\
{\leu0{ - (a + b*x)^(m + 1)*(c + d*x)^n *(e + f*x)^p = 0}}.

\tttdef\formx{exp(a + b*x)*(c + d*x)^m *(e + f*x)^n}
\section{Integrands involving \formx}

\subsection{Nondegenerate integrands}
\label{10A}

To exclude integrands with confluent roots, the following recurrences 
should be applied only if the resultant of the linear polynomials does 
not vanish: {\ttt{c*f - d*e \tneq 0}}.
\ads
{\len0{u*(c*f - d*e)*INT((A + B*x)}}\\
{\leu1{*exp(a + b*x)*(c + d*x)^m *(e + f*x)^n, x)}}\\
{\leu0{ - INT(((A*f - B*e)*u + V*(b*e + (n + 1)*f) + V*b*f*x)}}\\
{\leu1{*exp(a + b*x)*(c + d*x)^(m + 1)*(e + f*x)^n, x)}}\\
{\leu0{ + V*exp(a + b*x)*(c + d*x)^(m + 1)*(e + f*x)^(n + 1) = 0}},
\bds\bfl
where {\ttt{V = A*d - B*c}}, {\ttt{u = (m + 1)*d}}.
\efl\ads
{\len0{b*f*INT((A + B*x)*exp(a + b*x)*(c + d*x)^m *(e + f*x)^n, x)}}\\
{\leu0{ - INT((U*c - m*B*d*e + (U - m*B*f)*d*x)}}\\
{\leu1{*exp(a + b*x)*(c + d*x)^(m - 1)*(e + f*x)^n, x)}}\\
{\leu0{ - B*exp(a + b*x)*(c + d*x)^m *(e + f*x)^(n + 1) = 0}},
\bds\bfl
where {\ttt{U = (A*f - B*e)*b - (n + 1)*B*f}}.
\efl\ads
{\len0{u*d*INT((A + B*x)*exp(a + b*x)*(c + d*x)^m *(e + f*x)^n, x)}}\\
{\leu0{ - INT((B*u*e - V*(b*e + n*f) + (B*u - V*b)*f*x)}}\\
{\leu1{*exp(a + b*x)*(c + d*x)^(m + 1)*(e + f*x)^(n - 1), x)}}\\
{\leu0{ - V*exp(a + b*x)*(c + d*x)^(m + 1)*(e + f*x)^n = 0}},
\bds\bfl
where {\ttt{V = A*d - B*c}}, {\ttt{u = (m + 1)*d}}.
\efl

\subsection{Degenerate integrands}
\label{10B}

Now the roots of the two linear factors are required to coincide whence 
the corresponding resultant must vanish: {\ttt{c*f - d*e = 0}}. Note that 
each recurrence can be used backwards for its own (or its partner's) 
inverse.
\ads
{\len0{(m + n + 1)*d*INT(exp(a + b*x)*(c + d*x)^m *(e + f*x)^n, x)}}\\
{\leu0{ + b*INT(exp(a + b*x)*(c + d*x)^(m + 1)*(e + f*x)^n, x)}}\\
{\leu0{ - exp(a + b*x)*(c + d*x)^(m + 1)*(e + f*x)^n = 0}}.
\ads
{\len0{(m + n + 1)*d*INT(cos(a + b*x)*(c + d*x)^m *(e + f*x)^n, x)}}\\
{\leu0{ - b*INT(sin(a + b*x)*(c + d*x)^(m + 1)*(e + f*x)^n, x)}}\\
{\leu0{ - cos(a + b*x)*(c + d*x)^(m + 1)*(e + f*x)^n = 0}}.
\ads
{\len0{(m + n + 1)*d*INT(sin(a + b*x)*(c + d*x)^m *(e + f*x)^n, x)}}\\
{\leu0{ + b*INT(cos(a + b*x)*(c + d*x)^(m + 1)*(e + f*x)^n, x)}}\\
{\leu0{ - sin(a + b*x)*(c + d*x)^(m + 1)*(e + f*x)^n = 0}}.
\ads
{\len0{d*INT(exp(a + b*x)*(c + d*x)^m *(e + f*x)^n, x)}}\\
{\leu0{ - f*INT(exp(a + b*x)*(c + d*x)^(m + 1)*(e + f*x)^(n - 1), x) = 0}}.
\ads
{\len0{d*INT(cos(a + b*x)*(c + d*x)^m *(e + f*x)^n, x)}}\\
{\leu0{ - f*INT(cos(a + b*x)*(c + d*x)^(m + 1)*(e + f*x)^(n - 1), x) = 0}}.
\ads
{\len0{d*INT(sin(a + b*x)*(c + d*x)^m *(e + f*x)^n, x)}}\\
{\leu0{ - f*INT(sin(a + b*x)*(c + d*x)^(m + 1)*(e + f*x)^(n - 1), x) = 0}}.

\tttdef\formxi{(a + b*x)^m *(c + d*x)^n *(e + f*x + g*x^2)^p}
\section{Integrands involving \formxi}

Abbreviations used with these integrands:
\ads
{\reu0{ra = 2*a*g - b*f}}, {\ttt{rb = a*f - 2*b*e}}, 
{\ttt{rc = 2*c*g - d*f}}, {\ttt{rd = c*f - 2*d*e}}.\\[1ex]
{\reu0{re = 1/2 *(ra*d + rc*b)}}, {\ttt{rf = 1/2 *(ra*c + rd*b)}}, 
{\ttt{rg = 1/2 *(rb*c + rd*a)}}.\\[1ex]
{\reu0{se = 2*rf*g - re*f}}, {\ttt{sf = rg*g - re*e}}, 
{\ttt{sg = rg*f - 2*rf*e}}.

\subsection{Nondegenerate integrands}
\label{11A}

To exclude integrands with confluent roots, the following recurrences 
should be applied only if the overall discriminant does not vanish: 
{\ttt{(a*d - b*c)*(ra*a - rb*b)*(rc*c - rd*d)*(4*e*g - f^2) \tneq 0}}.
\ads
{\len0{u*1/2 *(ra*a - rb*b)*INT((A + B*x + C*x^2)}}\\
{\leu1{*(a + b*x)^m *(c + d*x)^n *(e + f*x + g*x^2)^p, x)}}\\
{\leu0{ + INT((((A*f - B*e)*b - U*a)*u - V*(v*e + (p + 1)*rd)}}\\
{\leu1{ + ((U*b - (B*g - C*f)*a)*u - V*(v*f + (p + 1)*rc))*x - V*v*g*x^2)}}\\
{\leu1{*(a + b*x)^(m + 1)*(c + d*x)^n *(e + f*x + g*x^2)^p, x)}}\\
{\leu0{ + V*(a + b*x)^(m + 1)*(c + d*x)^(n + 1)*(e + f*x + g*x^2)^(p + 1) = 0}},
\bds\bfl
where {\ttt{U = A*g - C*e}}, {\ttt{V = A*b^2 - B*a*b + C*a^2}}, 
{\ttt{u = (m + 1)*(a*d - b*c)}}, {\ttt{v = (m + n + 2*p + 4)*d}}.
\efl\ads
{\len0{t*g*INT((A + B*x + C*x^2)}}\\
{\leu1{*(a + b*x)^m *(c + d*x)^n *(e + f*x + g*x^2)^p, x)}}\\
{\leu0{ - INT((t*T*a + V*e - W*rg + (t*(T*b + U*a) + V*f - 2*W*rf)*x}}\\
{\leu1{ + (t*U*b + V*g - W*re)*x^2)}}\\
{\leu1{*(a + b*x)^(m - 1)*(c + d*x)^n *(e + f*x + g*x^2)^p, x)}}\\
{\leu0{ - C*(a + b*x)^m *(c + d*x)^(n + 1)*(e + f*x + g*x^2)^(p + 1) = 0}},
\bds\bfl
where {\ttt{T = A*g - C*e}}, {\ttt{U = B*g - C*f}}, 
{\ttt{V = (m + p + 1)*C*(a*d - b*c)}}, {\ttt{W = (p + 1)*C}}, 
{\ttt{t = (m + n + 2*p + 3)*d}}.
\efl\ads
{\len0{u*(4*e*g - f^2)*INT((A + B*x + C*x^2)}}\\
{\leu1{*(a + b*x)^m *(c + d*x)^n *(e + f*x + g*x^2)^p, x)}}\\
{\leu0{ - INT(((2*A*g - B*f + 2*C*e)*u + V*(v*b*c + w*a*d)}}\\
{\leu1{ + W*(a*c*g + (p + 1)*rf) + ((v + w)*V*b*d}}\\
{\leu1{ + W*(((v + 1)*b*c + (w + 1)*a*d)*g + (p + 1)*re))*x}}\\
{\leu1{ + (v + w + 1)*W*b*d*g*x^2)}}\\
{\leu1{*(a + b*x)^m *(c + d*x)^n *(e + f*x + g*x^2)^(p + 1), x)}}\\
{\leu0{ + (V + W*g*x)}}\\
{\leu1{*(a + b*x)^(m + 1)*(c + d*x)^(n + 1)*(e + f*x + g*x^2)^(p + 1) = 0}},
\bds\bfl
where {\ttt{V = (A*f - B*e)*se - (A*g - C*e)*sf}}, 
{\ttt{W = A*se - B*sf + C*sg}}, 
{\ttt{u = (p + 1)* 1/2 *(ra*a - rb*b)* 1/2 *(rc*c - rd*d)}}, 
{\ttt{v = m + p + 2}}, {\ttt{w = n + p + 2}}.
\efl\ads
{\len0{2*w*(w + b*d)*g*INT((A + B*x + C*x^2)}}\\
{\leu1{*(a + b*x)^m *(c + d*x)^n *(e + f*x + g*x^2)^p, x)}}\\
{\leu0{ + INT((V*(v*e + p*rg) - (U*e + W*a*c)*w}}\\
{\leu1{ + (V*(v*f + 2*p*rf) - (U*f + W*(a*d + b*c))*w)*x}}\\
{\leu1{ + (V*(v*g + p*re) - (U*g + W*b*d)*w)*x^2)}}\\
{\leu1{*(a + b*x)^m *(c + d*x)^n *(e + f*x + g*x^2)^(p - 1), x)}}\\
{\leu0{ - (V + C*w*f + 2*C*w*g*x)}}\\
{\leu1{*(a + b*x)^(m + 1)*(c + d*x)^(n + 1)*(e + f*x + g*x^2)^p = 0}},
\bds\bfl
where {\ttt{U = 2*(m + n + 2*p + 3)*(A*b*d - C*a*c)*g + %
C*((m + 1)*ra*c +{\fb} (n + 1)*rc*a)}}, {\ttt{V = 2*(m + n + 2*p + 3)*%
(B*b*d - C*(a*d + b*c))*g +{\fb} C*((m + 1)*ra*d + (n + 1)*rc*b)}}, 
{\ttt{W = p*C*(4*e*g - f^2)}}, 
{\ttt{v = (m + p + 1)*b*c + (n + p + 1)*a*d}}, 
{\ttt{w = (m + n + 2*p + 2)*b*d}}.
\efl\ads
{\len0{u*b*INT((A + B*x + C*x^2)}}\\
{\leu1{*(a + b*x)^m *(c + d*x)^n *(e + f*x + g*x^2)^p, x)}}\\
{\leu0{ - INT((U*c - V*(v*e + w*rd) + (U*d + C*u*c - V*(v*f + w*rc))*x}}\\
{\leu1{ + (C*u*d - V*v*g)*x^2)}}\\
{\leu1{*(a + b*x)^(m + 1)*(c + d*x)^(n - 1)*(e + f*x + g*x^2)^p, x)}}\\
{\leu0{ - 2*V*(a + b*x)^(m + 1)*(c + d*x)^n *(e + f*x + g*x^2)^(p + 1) = 0}},
\bds\bfl
where {\ttt{U = (m + 1)*(A*ra*b - (B*b - C*a)*rb)}}, 
{\ttt{V = A*b^2 - (B*b - C*a)*a}}, {\ttt{u = (m + 1)*(ra*a - rb*b)}}, 
{\ttt{v = 2*(m + n + 2*p + 3)*d}}, {\ttt{w = m + 2*p + 3}}.
\efl\ads
{\len0{u*1/2 *(rc*c - rd*d)*INT((A + B*x + C*x^2)}}\\
{\leu1{*(a + b*x)^m *(c + d*x)^n *(e + f*x + g*x^2)^p, x)}}\\
{\leu0{ - INT((U*u*a + V*(m*b*c + (n + 1)*a*d) + w*W*a*c + ((U*u + v*V*d)*b}}\\
{\leu1{ + W*((w + m)*b*c + (w + n + 1)*a*d))*x + (v + w)*W*b*d*x^2)}}\\
{\leu1{*(a + b*x)^(m - 1)*(c + d*x)^n *(e + f*x + g*x^2)^(p + 1), x)}}\\
{\leu0{ + (V + W*x)}}\\
{\leu1{*(a + b*x)^m *(c + d*x)^(n + 1)*(e + f*x + g*x^2)^(p + 1) = 0}},
\bds\bfl
where {\ttt{U = A*d^2 - B*c*d + C*c^2}}, 
{\ttt{V = (A*f - B*e)*rc - (A*g - C*e)*rd}}, 
{\ttt{W = (A*g - C*e)*rc - (B*g - C*f)*rd}}, 
{\ttt{u = (p + 1)*(4*e*g - f^2)}}, {\ttt{v = m + n + 1}}, 
{\ttt{w = 2*p + 3}}.
\efl\ads
{\len0{v*w*b^2*INT((A + B*x + C*x^2)}}\\
{\leu1{*(a + b*x)^m *(c + d*x)^n *(e + f*x + g*x^2)^p, x)}}\\
{\leu0{ - INT((U*e + V*v*(u*e + p*rd) - p*C*w*rg + (U*f + V*v*(u*f + p*rc)\\}}
{\leu1{ - 2*p*C*w*rf)*x + ((U + V*u*v)*g - p*C*w*re)*x^2)}}\\
{\leu1{*(a + b*x)^(m + 1)*(c + d*x)^n *(e + f*x + g*x^2)^(p - 1), x)}}\\
{\leu0{ + (V*v - C*w*a - C*w*b*x)}}\\
{\leu1{*(a + b*x)^(m + 1)*(c + d*x)^(n + 1)*(e + f*x + g*x^2)^p = 0}},
\bds\bfl
where {\ttt{U = (m + 1)*((m + n + 2*p + 3)*(A*d^2 - B*c*d + C*c^2)*b^2 -%
{\fb} (n + p + 1)*C*(a*d - b*c)^2)}}, {\ttt{V = A*b^2 - B*a*b + C*a^2}}, 
{\ttt{u = (n + 2*p + 1)*d}}, {\ttt{v = (m + n + 2*p + 3)*d}}, 
{\ttt{w = (m + 1)*(a*d - b*c)}}.
\efl

\subsection{Integrands with linear-linear cross-degeneracy}
\label{11B}

Here the roots of the two linear factors are required to coincide whence 
the corresponding resultant must vanish: {\ttt{a*d - b*c = 0}}. To 
exclude additional degeneracies, the recurrences should be applied only 
if also {\ttt{(ra*a - rb*b)*(4*e*g - f^2) \tneq 0}}.
\ads
{\len0{u*(ra*a - rb*b)*INT((A + B*x)}}\\
{\leu1{*(a + b*x)^m *(c + d*x)^n *(e + f*x + g*x^2)^p, x)}}\\
{\leu0{ - INT((u*(A*ra - B*rb) - v*V*f - 2*v*V*g*x)}}\\
{\leu1{*(a + b*x)^(m + 1)*(c + d*x)^n *(e + f*x + g*x^2)^p, x)}}\\
{\leu0{ - 2*V*(a + b*x)^(m + 1)*(c + d*x)^n *(e + f*x + g*x^2)^(p + 1) = 0}},
\bds\bfl
where {\ttt{V = A*b - B*a}}, {\ttt{u = m + n + 1}}, 
{\ttt{v = m + n + 2*p + 3}}.
\efl\ads
{\len0{2*u*g*INT((A + B*x)}}\\
{\leu1{*(a + b*x)^m *(c + d*x)^n *(e + f*x + g*x^2)^p, x)}}\\
{\leu0{ - INT((u*U*a + V*rb + (u*U*b + V*ra)*x)}}\\
{\leu1{*(a + b*x)^(m - 1)*(c + d*x)^n *(e + f*x + g*x^2)^p, x)}}\\
{\leu0{ - 2*B*(a + b*x)^m *(c + d*x)^n *(e + f*x + g*x^2)^(p + 1) = 0}},
\bds\bfl
where {\ttt{U = 2*A*g - B*f}}, {\ttt{V = (m + n)*B}}, 
{\ttt{u = m + n + 2*p + 2}}.
\efl\ads
{\len0{(p + 1)* 1/2 *(ra*a - rb*b)*(4*e*g - f^2)*INT((A + B*x)}}\\
{\leu1{*(a + b*x)^m *(c + d*x)^n *(e + f*x + g*x^2)^p, x)}}\\
{\leu0{ + INT((U*u - V*((p + 1)*ra + a*g) - V*(u + b)*g*x)}}\\
{\leu1{*(a + b*x)^m *(c + d*x)^n *(e + f*x + g*x^2)^(p + 1), x)}}\\
{\leu0{ - (U - V*g*x)}}\\
{\leu1{*(a + b*x)^(m + 1)*(c + d*x)^n *(e + f*x + g*x^2)^(p + 1) = 0}},
\bds\bfl
where {\ttt{U = A*rb*g - (A*f - B*e)*ra}}, {\ttt{V = A*ra - B*rb}}, 
{\ttt{u = (m + n + 2*p + 3)*b}}.
\efl\ads
{\len0{2*v*(v + b)*g*INT((A + B*x)}}\\
{\leu1{*(a + b*x)^m *(c + d*x)^n *(e + f*x + g*x^2)^p, x)}}\\
{\leu0{ + p*INT((U*rb - V*v*a + (U*ra - V*v*b)*x)}}\\
{\leu1{*(a + b*x)^m *(c + d*x)^n *(e + f*x + g*x^2)^(p - 1), x)}}\\
{\leu0{ - (U + B*v*f + 2*B*v*g*x)}}\\
{\leu1{*(a + b*x)^(m + 1)*(c + d*x)^n *(e + f*x + g*x^2)^p = 0}},
\bds\bfl
where {\ttt{U = 2*(m + n + 2*p + 2)*(A*b - B*a)*g + (m + n + 1)*B*ra}}, 
{\ttt{V = B*(4*e*g - f^2)}}, {\ttt{v = (m + n + 2*p + 1)*b}}.
\efl\ads
{\len0{b*INT((A + B*x)*(a + b*x)^m *(c + d*x)^n *(e + f*x + g*x^2)^p, x)}}\\
{\leu0{ - d*INT((A + B*x)}}\\
{\leu1{*(a + b*x)^(m + 1)*(c + d*x)^(n - 1)*(e + f*x + g*x^2)^p, x) = 0}}.
\ads
{\len0{(p + 1)*(4*e*g - f^2)*INT((A + B*x)}}\\
{\leu1{*(a + b*x)^m *(c + d*x)^n *(e + f*x + g*x^2)^p, x)}}\\
{\leu0{ - INT(((m + n)*U*b + v*V*a + (v + m + n)*V*b*x)}}\\
{\leu1{*(a + b*x)^(m - 1)*(c + d*x)^n *(e + f*x + g*x^2)^(p + 1), x)}}\\
{\leu0{ + (U + V*x)*(a + b*x)^m *(c + d*x)^n *(e + f*x + g*x^2)^(p + 1) = 0}},
\bds\bfl
where {\ttt{U = A*f - 2*B*e}}, {\ttt{V = 2*A*g - B*f}}, 
{\ttt{v = 2*p + 3}}.
\efl\ads
{\len0{u*v*b^2*INT((A + B*x)}}\\
{\leu1{*(a + b*x)^m *(c + d*x)^n *(e + f*x + g*x^2)^p, x)}}\\
{\leu0{ + p*INT((u*U*f + v*B*rb + (2*u*U*g + v*B*ra)*x)}}\\
{\leu1{*(a + b*x)^(m + 1)*(c + d*x)^n *(e + f*x + g*x^2)^(p - 1), x)}}\\
{\leu0{ - (u*U + v*B*a + v*B*b*x)}}\\
{\leu1{*(a + b*x)^(m + 1)*(c + d*x)^n *(e + f*x + g*x^2)^p = 0}},
\bds\bfl
where {\ttt{U = A*b - B*a}}, {\ttt{u = m + n + 2*p + 2}}, 
{\ttt{v = m + n + 1}}.
\efl

\subsection{Integrands with linear-quadratic cross-degeneracy}
\label{11C}

\begin{sloppypar}
Here the root of the first linear factor is required to coincide with one 
root of the quadratic whence the corresponding resultant must vanish: 
{\ttt{ra*a - rb*b = 0}}. To exclude additional degeneracies, the 
recurrences should be applied only if also {\ttt{(a*d - b*c)*{\fb}%
(4*e*g - f^2) \tneq 0}}.
\end{sloppypar}
\ads
{\len0{u*(a*d - b*c)*INT((A + B*x)}}\\
{\leu1{*(a + b*x)^m *(c + d*x)^n *(e + f*x + g*x^2)^p, x)}}\\
{\leu0{ - INT(((A*d - B*c)*u + (n + 1)*(A*re - B*rf) - V*v*c - V*v*d*x)}}\\
{\leu1{*(a + b*x)^(m + 1)*(c + d*x)^n *(e + f*x + g*x^2)^p, x)}}\\
{\leu0{ - V*b*(a + b*x)^m *(c + d*x)^(n + 1)*(e + f*x + g*x^2)^(p + 1) = 0}},
\bds\bfl
where {\ttt{V = A*b - B*a}}, {\ttt{u = (m + p + 1)*ra}}, 
{\ttt{v = (m + n + 2*p + 3)*g}}.
\efl\ads
{\len0{u*d*g*INT((A + B*x)*(a + b*x)^m *(c + d*x)^n *(e + f*x + g*x^2)^p, x)}}\\
{\leu0{ - INT((u*U*a + V*c + W*rf + (u*U*b + V*d + W*re)*x)}}\\
{\leu1{*(a + b*x)^(m - 1)*(c + d*x)^n *(e + f*x + g*x^2)^p, x)}}\\
{\leu0{ - B*b}}\\
{\leu1{*(a + b*x)^(m - 1)*(c + d*x)^(n + 1)*(e + f*x + g*x^2)^(p + 1) = 0}},
\bds\bfl
where {\ttt{U = (A*d - B*c)*g}}, {\ttt{V = (m + p)*B*ra}}, 
{\ttt{W = (n + 1)*B}}, {\ttt{u = m + n + 2*p + 2}}.
\efl\ads
{\len0{(n + 1)* 1/2 *(rc*c - rd*d)*b*INT((A + B*x)}}\\
{\leu1{*(a + b*x)^m *(c + d*x)^n *(e + f*x + g*x^2)^p, x)}}\\
{\leu0{ - INT(((n + 1)*(A*re - B*rf) - V*(v*a - (m + p + 1)*ra) - V*v*b*x)}}\\
{\leu1{*(a + b*x)^m *(c + d*x)^(n + 1)*(e + f*x + g*x^2)^p, x)}}\\
{\leu0{ - V*b*(a + b*x)^m *(c + d*x)^(n + 1)*(e + f*x + g*x^2)^(p + 1) = 0}},
\bds\bfl
where {\ttt{V = A*d - B*c}}, {\ttt{v = (m + n + 2*p + 3)*g}}.
\efl\ads
{\len0{(m + n + 2*p + 2)*b*g*INT((A + B*x)}}\\
{\leu1{*(a + b*x)^m *(c + d*x)^n *(e + f*x + g*x^2)^p, x)}}\\
{\leu0{ - INT((U*c + n*B*rf + (U*d + n*B*re)*x)}}\\
{\leu1{*(a + b*x)^m *(c + d*x)^(n - 1)*(e + f*x + g*x^2)^p, x)}}\\
{\leu0{ - B*b*(a + b*x)^m *(c + d*x)^n *(e + f*x + g*x^2)^(p + 1) = 0}},
\bds\bfl
where {\ttt{U = (m + n + 2*p + 2)*(A*b - B*a)*g + (m + p + 1)*B*ra}}.
\efl\ads
{\len0{(p + 1)*u*(rc*c - rd*d)*INT((A + B*x)}}\\
{\leu1{*(a + b*x)^m *(c + d*x)^n *(e + f*x + g*x^2)^p, x)}}\\
{\leu0{ - INT((((w + n)*U*u + (n + 1)*(V - W*f))*d - 2*(w + m)*W*c*g}}\\
{\leu1{ - 2*(w + m + n + 1)*W*d*g*x)}}\\
{\leu1{*(a + b*x)^m *(c + d*x)^n *(e + f*x + g*x^2)^(p + 1), x)}}\\
{\leu0{ + (U*u + V - W*f - 2*W*g*x)}}\\
{\leu1{*(a + b*x)^m *(c + d*x)^(n + 1)*(e + f*x + g*x^2)^(p + 1) = 0}},
\bds\bfl
where {\ttt{U = A*d - B*c}}, {\ttt{V = m*(A*se - B*sf)}}, 
{\ttt{W = m*(A*re - B*rf) -{\fb} (m + p + 1)*(A*rc - B*rd)*b}}, 
{\ttt{u = (m + p + 1)*(4*e*g - f^2)*b}}, {\ttt{w = 2*p + 3}}.
\efl\ads
{\len0{2*v*(v + d)*b*g*INT((A + B*x)}}\\
{\leu1{*(a + b*x)^m *(c + d*x)^n *(e + f*x + g*x^2)^p, x)}}\\
{\leu0{ - INT((U*(m*rf - w*rd) + V*sf + W*w*c}}\\
{\leu1{ + (U*(m*re - w*rc) + V*se + W*w*d)*x)}}\\
{\leu1{*(a + b*x)^m *(c + d*x)^n *(e + f*x + g*x^2)^(p - 1), x)}}\\
{\leu0{ - (U*b + B*(m*ra*d + v*b*f) + 2*B*v*b*g*x)}}\\
{\leu1{*(a + b*x)^m *(c + d*x)^(n + 1)*(e + f*x + g*x^2)^p = 0}},
\bds\bfl
where {\ttt{U = 2*(m + n + 2*p + 2)*(A*d - B*c)*g + (n + 1)*B*rc}}, 
{\ttt{V = m*(n + 1)*B*d}}, {\ttt{W = (n + 2*p + 1)*B*(4*e*g - f^2)*d}}, 
{\ttt{v = (m + n + 2*p + 1)*d}}, {\ttt{w = (m + p)*b}}.
\efl\ads
{\len0{(m + p + 1)*ra*b*INT((A + B*x)}}\\
{\leu1{*(a + b*x)^m *(c + d*x)^n *(e + f*x + g*x^2)^p, x)}}\\
{\leu0{ + INT((n*(A*re - B*rf) - V*c - V*d*x)}}\\
{\leu1{*(a + b*x)^(m + 1)*(c + d*x)^(n - 1)*(e + f*x + g*x^2)^p, x)}}\\
{\leu0{ + (A*b - B*a)*b}}\\
{\leu1{*(a + b*x)^m *(c + d*x)^n *(e + f*x + g*x^2)^(p + 1) = 0}},
\bds\bfl
where {\ttt{V = (m + n + 2*p + 2)*(A*b - B*a)*g + (m + p + 1)*B*ra}}.
\efl\ads
{\len0{(n + 1)*re*d*INT((A + B*x)}}\\
{\leu1{*(a + b*x)^m *(c + d*x)^n *(e + f*x + g*x^2)^p, x)}}\\
{\leu0{ - b*INT((U*rf + V*(v*a - (m + p)*ra) + (U*re + V*v*b)*x)}}\\
{\leu1{*(a + b*x)^(m - 1)*(c + d*x)^(n + 1)*(e + f*x + g*x^2)^p, x)}}\\
{\leu0{ + V*b^2}}\\
{\leu1{*(a + b*x)^(m - 1)*(c + d*x)^(n + 1)*(e + f*x + g*x^2)^(p + 1) = 0}},
\bds\bfl
where {\ttt{U = (n + 1)*B}}, {\ttt{V = A*d - B*c}}, 
{\ttt{v = (m + n + 2*p + 2)*g}}.
\efl\ads
{\len0{(p + 1)*se*INT((A + B*x)}}\\
{\leu1{*(a + b*x)^m *(c + d*x)^n *(e + f*x + g*x^2)^p, x)}}\\
{\leu0{ - b*INT(((n + 1)*(A*re - B*rf) - (n + p + 2)*(A*d - B*c)*ra}}\\
{\leu1{ - v*V*c - v*V*d*x)\\}}
{\leu1{*(a + b*x)^(m - 1)*(c + d*x)^n *(e + f*x + g*x^2)^(p + 1), x)}}\\
{\leu0{ - V*(a + b*x)^m *(c + d*x)^(n + 1)*(e + f*x + g*x^2)^(p + 1) = 0}},
\bds\bfl
where {\ttt{V = A*b*g + B*(a*g - b*f)}}, {\ttt{v = m + n + 2*p + 3}}.
\efl\ads
{\len0{u*b^2*d*INT((A + B*x)}}\\
{\leu1{*(a + b*x)^m *(c + d*x)^n *(e + f*x + g*x^2)^p, x)}}\\
{\leu0{ + INT((u*U*(a*g - b*f) + V*c - W*rf - (u*U*b*g - V*d + W*re)*x)}}\\
{\leu1{*(a + b*x)^(m + 1)*(c + d*x)^n *(e + f*x + g*x^2)^(p - 1), x)}}\\
{\leu0{ - B*b*(a + b*x)^(m + 1)*(c + d*x)^(n + 1)*(e + f*x + g*x^2)^p = 0}},
\bds\bfl
where {\ttt{U = A*d - B*c}}, {\ttt{V = (n + p + 1)*B*ra}}, 
{\ttt{W = (n + 1)*B}}, {\ttt{u = m + n + 2*p + 2}}.
\efl\ads
{\len0{2*(m + p + 1)*(p + 1)*(4*e*g - f^2)*b*INT((A + B*x)}}\\
{\leu1{*(a + b*x)^m *(c + d*x)^n *(e + f*x + g*x^2)^p, x)}}\\
{\leu0{ - INT((n*U*d + v*V*c + (v + n)*V*d*x)}}\\
{\leu1{*(a + b*x)^m *(c + d*x)^(n - 1)*(e + f*x + g*x^2)^(p + 1), x)}}\\
{\leu0{ + (U + V*x)*(a + b*x)^m *(c + d*x)^n *(e + f*x + g*x^2)^(p + 1) = 0}},
\bds\bfl
where {\ttt{U = m*A*ra + (m + 2*p + 2)*(A*f - 2*B*e)*b}}, 
{\ttt{V = (m + 2*p + 2)*(2*A*g - B*f)*b + m*B*ra}}, 
{\ttt{v = m + 2*p + 3}}.
\efl\ads
{\len0{(m + n + 2*p + 2)*v*b*d*INT((A + B*x)}}\\
{\leu1{*(a + b*x)^m *(c + d*x)^n *(e + f*x + g*x^2)^p, x)}}\\
{\leu0{ + INT((U*(u*f - w) - B*u*v*e + (U*u*g - B*v*w)*x)}}\\
{\leu1{*(a + b*x)^m *(c + d*x)^(n + 1)*(e + f*x + g*x^2)^(p - 1), x)}}\\
{\leu0{ - b*(U + B*v*x)}}\\
{\leu1{*(a + b*x)^m *(c + d*x)^(n + 1)*(e + f*x + g*x^2)^p = 0}},
\bds\bfl
where {\ttt{U = (m + n + 2*p + 2)*A*d - (m + 2*p + 1)*B*c}}, 
{\ttt{u = (m + 2*p)*b}}, {\ttt{v = (n + 1)*d}}, 
{\ttt{w = m*a*g + p*b*f}}.
\efl

\subsection{Integrands with a self-degenerate quadratic}
\label{11D}

Here the quadratic polynomial is required to possess a double root whence 
its discriminant must vanish: {\ttt{4*e*g - f^2 = 0}}. To exclude 
additional degeneracies, the recurrences should be applied only if also 
{\ttt{ra*rc \tneq 0}}.
\ads
{\len0{u*(a*d - b*c)*INT((A + B*x)}}\\
{\leu1{*(a + b*x)^m *(c + d*x)^n *(e + f*x + g*x^2)^p, x)}}\\
{\leu0{ - INT(((A*d - B*c)*u - V*(v*c - (n + 1)*rc) - V*v*d*x)}}\\
{\leu1{*(a + b*x)^(m + 1)*(c + d*x)^n *(e + f*x + g*x^2)^p, x)}}\\
{\leu0{ - V*(f + 2*g*x)}}\\
{\leu1{*(a + b*x)^(m + 1)*(c + d*x)^(n + 1)*(e + f*x + g*x^2)^p = 0}},
\bds\bfl
where {\ttt{V = A*b - B*a}}, {\ttt{u = (m + 1)*ra}}, 
{\ttt{v = 2*(m + n + 2*p + 3)*g}}.
\efl\ads
{\len0{2*(m + n + 2*p + 2)*d*g*INT((A + B*x)}}\\
{\leu1{*(a + b*x)^m *(c + d*x)^n *(e + f*x + g*x^2)^p, x)}}\\
{\leu0{ - INT((U*a + V*f + (U*b + 2*V*g)*x)}}\\
{\leu1{*(a + b*x)^(m - 1)*(c + d*x)^n *(e + f*x + g*x^2)^p, x)}}\\
{\leu0{ - B*(f + 2*g*x)*(a + b*x)^m *(c + d*x)^(n + 1)*(e + f*x + g*x^2)^p = 0}},
\bds\bfl
where {\ttt{U = 2*(m + n + 2*p + 2)*(A*d - B*c)*g + (m + n + 1)*B*rc}}, 
{\ttt{V = m*B*(a*d - b*c)}}.
\efl\ads
{\len0{u*w*INT((A + B*x)*(a + b*x)^m *(c + d*x)^n *(e + f*x + g*x^2)^p, x)}}\\
{\leu0{ + INT((v*(U*a*d + V*b*c - A*w*b*d)}}\\
{\leu1{ - (U + V - B*w)*((v + m + 1)*b*c + (v + n + 1)*a*d)}}\\
{\leu1{ - (v + m + n + 2)*(U + V - B*w)*b*d*x)}}\\
{\leu1{*(a + b*x)^m *(c + d*x)^n *(e + f*x + g*x^2)^(p + 1), x)}}\\
{\leu0{ - u*(A*f - 2*B*e + (2*A*g - B*f)*x)}}\\
{\leu1{*(a + b*x)^(m + 1)*(c + d*x)^(n + 1)*(e + f*x + g*x^2)^p}}\\
{\leu0{ + (U + V - B*w)}}\\
{\leu1{*(a + b*x)^(m + 1)*(c + d*x)^(n + 1)*(e + f*x + g*x^2)^(p + 1) = 0}},
\bds\bfl
where {\ttt{U = (m + 2*p + 2)*(A*rc - B*rd)*b}}, 
{\ttt{V = (n + 2*p + 2)*(A*ra - B*rb)*d}}, 
{\ttt{u = (p + 1)*(ra*c - rb*d)}}, {\ttt{v = 2*(p + 1)}}, 
{\ttt{w = (2*p + 1)*(ra*c - rb*d)}}.
\efl\ads
{\len0{2*t*(t + b*d)*INT((A + B*x)}}\\
{\leu1{*(a + b*x)^m *(c + d*x)^n *(e + f*x + g*x^2)^p, x)}}\\
{\leu0{ + INT((V*u + v*(V*b + W*a)*rd + w*(V*d + W*c)*rb}}\\
{\leu1{ + (v*(V*b + W*a)*rc + w*(V*d + W*c)*ra - W*u)*x)}}\\
{\leu1{*(a + b*x)^m *(c + d*x)^n *(e + f*x + g*x^2)^(p - 1), x)}}\\
{\leu0{ - (V*f + 2*W*e + (2*V*g + W*f)*x)}}\\
{\leu1{*(a + b*x)^(m + 1)*(c + d*x)^(n + 1)*(e + f*x + g*x^2)^(p - 1)}}\\
{\leu0{ - 2*B*t*(a + b*x)^(m + 1)*(c + d*x)^(n + 1)*(e + f*x + g*x^2)^p = 0}},
\bds\bfl
where {\ttt{V = (m + n + 2*p + 2)*A*b*d - B*((m + 2*p + 1)*b*c + %
(n + 2*p + 1)*a*d)}}, {\ttt{W = 2*p*B*b*d}}, 
{\ttt{t = (m + n + 2*p + 1)*b*d}}, {\ttt{u = (2*p - 1)*(ra*c - rb*d)}}, 
{\ttt{v = m + 2*p}}, {\ttt{w = n + 2*p}}.
\efl\ads
{\len0{u*b*INT((A + B*x)*(a + b*x)^m *(c + d*x)^n *(e + f*x + g*x^2)^p, x)}}\\
{\leu0{ - INT(((B*u + V*v)*c - n*V*rc + (B*u + V*v)*d*x)}}\\
{\leu1{*(a + b*x)^(m + 1)*(c + d*x)^(n - 1)*(e + f*x + g*x^2)^p, x)}}\\
{\leu0{ + V*(f + 2*g*x)*(a + b*x)^(m + 1)*(c + d*x)^n *(e + f*x + g*x^2)^p = 0}},
\bds\bfl
where {\ttt{V = A*b - B*a}}, {\ttt{u = (m + 1)*ra}}, 
{\ttt{v = 2*(m + n + 2*p + 2)*g}}.
\efl\ads
{\len0{(2*p + 1)*u*rc*INT((A + B*x)}}\\
{\leu1{*(a + b*x)^m *(c + d*x)^n *(e + f*x + g*x^2)^p, x)}}\\
{\leu0{ - INT((V*(v*a*d - m*w) + W*(v*b*c + (n + 1)*w) + v*(V + W)*b*d*x)}}\\
{\leu1{*(a + b*x)^(m - 1)*(c + d*x)^n *(e + f*x + g*x^2)^(p + 1), x)}}\\
{\leu0{ - u*(A*f - 2*B*e + (2*A*g - B*f)*x)}}\\
{\leu1{*(a + b*x)^m *(c + d*x)^(n + 1)*(e + f*x + g*x^2)^p}}\\
{\leu0{ + (V + W)*(a + b*x)^m *(c + d*x)^(n + 1)*(e + f*x + g*x^2)^(p + 1) = 0}},
\bds\bfl
where {\ttt{V = (n + 1)*(A*ra - B*rb)*d + (2*p + 1)*(A*d - B*c)*ra}}, 
{\ttt{W = m*(A*rc - B*rd)*b}}, {\ttt{u = (p + 1)*(ra*c - rb*d)}}, 
{\ttt{v = m + n + 2*p + 3}}, {\ttt{w = a*d - b*c}}.
\efl\ads
{\len0{(m + n + 2*p + 2)*v*w*b*d*INT((A + B*x)}}\\
{\leu1{*(a + b*x)^m *(c + d*x)^n *(e + f*x + g*x^2)^p, x)}}\\
{\leu0{ + INT((U*((n + 1)*rb*d + (2*p - 1)*ra*c) + (U + p*B*v)*w*rd}}\\
{\leu1{ + (U*(n + 2*p)*ra*d + (U + p*B*v)*w*rc)*x)}}\\
{\leu1{*(a + b*x)^(m + 1)*(c + d*x)^n *(e + f*x + g*x^2)^(p - 1), x)}}\\
{\leu0{ - U*(rb + ra*x)}}\\
{\leu1{*(a + b*x)^(m + 1)*(c + d*x)^(n + 1)*(e + f*x + g*x^2)^(p - 1)}}\\
{\leu0{ - B*v*w*(a + b*x)^(m + 1)*(c + d*x)^(n + 1)*(e + f*x + g*x^2)^p = 0}},
\bds\bfl
where {\ttt{U = (m + 1)*(A*d - B*c)*b + (n + 2*p + 1)*(A*b - B*a)*d}}, 
{\ttt{v = 2*(a*d - b*c)}}, {\ttt{w = (m + 1)*b}}.
\efl

\subsection{Doubly degenerate integrands}
\label{11E}

Now the roots of both linear factors and one root of the quadratic are 
required to coincide whence the following must vanish: 
{\ttt{a*d - b*c = ra*a - rb*b = 0}}. To exclude additional degeneracies, 
the recurrences should be applied only if also {\ttt{ra \tneq 0}}. Note 
that each relation can be used backwards for its own inverse.
\ads
{\len0{u*ra*INT((a + b*x)^m *(c + d*x)^n *(e + f*x + g*x^2)^p, x)}}\\
{\leu0{ - (u + p + 1)*g}}\\
{\leu1{*INT((a + b*x)^(m + 1)*(c + d*x)^n *(e + f*x + g*x^2)^p, x)}}\\
{\leu0{ + b*(a + b*x)^m *(c + d*x)^n *(e + f*x + g*x^2)^(p + 1) = 0}},
\bds\bfl
where {\ttt{u = m + n + p + 1}}.
\efl\ads
{\len0{(m + n + p + 1)*(p + 1)*(4*e*g - f^2)*b}}\\
{\leu1{*INT((a + b*x)^m *(c + d*x)^n *(e + f*x + g*x^2)^p, x)}}\\
{\leu0{ - w*(w + 1)*b*g}}\\
{\leu1{*INT((a + b*x)^m *(c + d*x)^n *(e + f*x + g*x^2)^(p + 1), x)}}\\
{\leu0{ - ((p + 1)*ra - w*a*g - w*b*g*x)}}\\
{\leu1{*(a + b*x)^m *(c + d*x)^n *(e + f*x + g*x^2)^(p + 1) = 0}},
\bds\bfl
where {\ttt{w = m + n + 2*p + 2}}.
\efl\ads
{\len0{b*INT((a + b*x)^m *(c + d*x)^n *(e + f*x + g*x^2)^p, x)}}\\
{\leu0{ - d*INT((a + b*x)^(m + 1)*(c + d*x)^(n - 1)*(e + f*x + g*x^2)^p, x)}}\\
{\leu0{ = 0}}.
\ads
{\len0{(p + 1)*ra*INT((a + b*x)^m *(c + d*x)^n *(e + f*x + g*x^2)^p, x)}}\\
{\leu0{ + (m + n + 2*p + 2)*b^2}}\\
{\leu1{*INT((a + b*x)^(m - 1)*(c + d*x)^n *(e + f*x + g*x^2)^(p + 1), x)}}\\
{\leu0{ - b*(a + b*x)^m *(c + d*x)^n *(e + f*x + g*x^2)^(p + 1) = 0}}.
\bds

Now the roots of the two linear factors are required to coincide while 
the quadratic is required to possess a double root whence the following{\pagebreak[1]} 
must vanish: {\ttt{a*d - b*c = 4*e*g - f^2 = 0}}. To exclude additional 
degeneracies, the recurrences should be applied only if also 
{\ttt{ra \tneq 0}}. Note that each relation can be used backwards for its 
own inverse.
\ads
{\len0{(m + n + 1)*ra*INT((a + b*x)^m *(c + d*x)^n *(e + f*x + g*x^2)^p, x)}}\\
{\leu0{ - 2*(m + n + 2*p + 2)*g}}\\
{\leu1{*INT((a + b*x)^(m + 1)*(c + d*x)^n *(e + f*x + g*x^2)^p, x)}}\\
{\leu0{ + (f + 2*g*x)*(a + b*x)^(m + 1)*(c + d*x)^n *(e + f*x + g*x^2)^p = 0}}.
\ads
{\len0{(p + 1)*(2*p + 1)*(ra*a - rb*b)}}\\
{\leu1{*INT((a + b*x)^m *(c + d*x)^n *(e + f*x + g*x^2)^p, x)}}\\
{\leu0{ - w*(w + b)}}\\
{\leu1{*INT((a + b*x)^m *(c + d*x)^n *(e + f*x + g*x^2)^(p + 1), x)}}\\
{\leu0{ - (p + 1)*(rb + ra*x)}}\\
{\leu1{*(a + b*x)^(m + 1)*(c + d*x)^n *(e + f*x + g*x^2)^p}}\\
{\leu0{ + w*(a + b*x)^(m + 1)*(c + d*x)^n *(e + f*x + g*x^2)^(p + 1) = 0}},
\bds\bfl
where {\ttt{w = (m + n + 2*p + 2)*b}}.
\efl\ads
{\len0{b*INT((a + b*x)^m *(c + d*x)^n *(e + f*x + g*x^2)^p, x)}}\\
{\leu0{ - d*INT((a + b*x)^(m + 1)*(c + d*x)^(n - 1)*(e + f*x + g*x^2)^p, x)}}\\
{\leu0{ = 0}}.
\ads
{\len0{(p + 1)*(2*p + 1)*ra}}\\
{\leu1{*INT((a + b*x)^m *(c + d*x)^n *(e + f*x + g*x^2)^p, x)}}\\
{\leu0{ - (m + n + 2*p + 2)*w*b}}\\
{\leu1{*INT((a + b*x)^(m - 1)*(c + d*x)^n *(e + f*x + g*x^2)^(p + 1), x)}}\\
{\leu0{ - (p + 1)*(rb + ra*x)*(a + b*x)^m *(c + d*x)^n *(e + f*x + g*x^2)^p}}\\
{\leu0{ + w*(a + b*x)^m *(c + d*x)^n *(e + f*x + g*x^2)^(p + 1) = 0}},
\bds\bfl
where {\ttt{w = (m + n)*b}}.
\efl

Now the roots of the two linear factors are required to coincide with 
different roots of the quadratic whence the following must vanish: 
{\ttt{re = rf = 0}}. To exclude additional degeneracies, the recurrences 
should be applied only if also {\ttt{a*d - b*c \tneq 0}}. Note that each 
relation can be used backwards for its own inverse.
\ads
{\len0{(m + p + 1)*(a*d - b*c)}}\\
{\leu1{*INT((a + b*x)^m *(c + d*x)^n *(e + f*x + g*x^2)^p, x)}}\\
{\leu0{ - (m + n + 2*p + 2)*d}}\\
{\leu1{*INT((a + b*x)^(m + 1)*(c + d*x)^n *(e + f*x + g*x^2)^p, x)}}\\
{\leu0{ + (a + b*x)^(m + 1)*(c + d*x)^(n + 1)*(e + f*x + g*x^2)^p = 0}}.
\ads
{\len0{u*v*(ra*c - rb*d)}}\\
{\leu1{*INT((a + b*x)^m *(c + d*x)^n *(e + f*x + g*x^2)^p, x)}}\\
{\leu0{ - (u + v)*(u + v + 1)*b*d}}\\
{\leu1{*INT((a + b*x)^m *(c + d*x)^n *(e + f*x + g*x^2)^(p + 1), x)}}\\
{\leu0{ + (u*a*d + v*b*c + (u + v)*b*d*x)}}\\
{\leu1{*(a + b*x)^m *(c + d*x)^n *(e + f*x + g*x^2)^(p + 1) = 0}},
\bds\bfl
where {\ttt{u = m + p + 1}}, {\ttt{v = n + p + 1}}.
\efl\ads
{\len0{(m + p + 1)*b*INT((a + b*x)^m *(c + d*x)^n *(e + f*x + g*x^2)^p, x)}}\\
{\leu0{ + (n + p)*d}}\\
{\leu1{*INT((a + b*x)^(m + 1)*(c + d*x)^(n - 1)*(e + f*x + g*x^2)^p, x)}}\\
{\leu0{ - (a + b*x)^(m + 1)*(c + d*x)^n *(e + f*x + g*x^2)^p = 0}}.
\ads
{\len0{(n + p + 1)*ra*INT((a + b*x)^m *(c + d*x)^n *(e + f*x + g*x^2)^p, x)}}\\
{\leu0{ + (m + n + 2*p + 2)*b^2}}\\
{\leu1{*INT((a + b*x)^(m - 1)*(c + d*x)^n *(e + f*x + g*x^2)^(p + 1), x)}}\\
{\leu0{ - b*(a + b*x)^m *(c + d*x)^n *(e + f*x + g*x^2)^(p + 1) = 0}}.
\bds

Now the quadratic is required to possess a double root which coincides 
with the root of the first linear factor whence the following must 
vanish: {\ttt{ra = rb = 0}}. To exclude additional degeneracies, the 
recurrences should be applied only if also {\ttt{a*d - b*c \tneq 0}}. 
Note that each relation can be used backwards for its own inverse.
\ads
{\len0{u*(a*d - b*c)*INT((a + b*x)^m *(c + d*x)^n *(e + f*x + g*x^2)^p, x)}}\\
{\leu0{ - (u + n + 1)*d}}\\
{\leu1{*INT((a + b*x)^(m + 1)*(c + d*x)^n *(e + f*x + g*x^2)^p, x)}}\\
{\leu0{ + (a + b*x)^(m + 1)*(c + d*x)^(n + 1)*(e + f*x + g*x^2)^p = 0}},
\bds\bfl
where {\ttt{u = m + 2*p + 1}}.
\efl\ads
{\len0{(n + 1)*(a*d - b*c)}}\\
{\leu1{*INT((a + b*x)^m *(c + d*x)^n *(e + f*x + g*x^2)^p, x)}}\\
{\leu0{ + (m + n + 2*p + 2)*b}}\\
{\leu1{*INT((a + b*x)^m *(c + d*x)^(n + 1)*(e + f*x + g*x^2)^p, x)}}\\
{\leu0{ - (a + b*x)^(m + 1)*(c + d*x)^(n + 1)*(e + f*x + g*x^2)^p = 0}}.
\ads
{\len0{(m + 2*p + 1)*u*1/2 *(rc*c - rd*d)}}\\
{\leu1{*INT((a + b*x)^m *(c + d*x)^n *(e + f*x + g*x^2)^p, x)}}\\
{\leu0{ - v*(v + d)}}\\
{\leu1{*INT((a + b*x)^m *(c + d*x)^n *(e + f*x + g*x^2)^(p + 1), x)}}\\
{\leu0{ + (u*(a*d - b*c) + v*a + v*b*x)}}\\
{\leu1{*(a + b*x)^(m - 1)*(c + d*x)^(n + 1)*(e + f*x + g*x^2)^(p + 1) = 0}},
\bds\bfl
where {\ttt{u = m + 2*p + 2}}, {\ttt{v = (m + n + 2*p + 2)*d}}.
\efl\ads
{\len0{(m + 2*p + 1)*b*INT((a + b*x)^m *(c + d*x)^n *(e + f*x + g*x^2)^p, x)}}\\
{\leu0{ + n*d}}\\
{\leu1{*INT((a + b*x)^(m + 1)*(c + d*x)^(n - 1)*(e + f*x + g*x^2)^p, x)}}\\
{\leu0{ - (a + b*x)^(m + 1)*(c + d*x)^n *(e + f*x + g*x^2)^p = 0}}.
\ads
{\len0{u*rc*INT((a + b*x)^m *(c + d*x)^n *(e + f*x + g*x^2)^p, x)}}\\
{\leu0{ + 2*(u + n + 1)*b*d}}\\
{\leu1{*INT((a + b*x)^(m - 1)*(c + d*x)^n *(e + f*x + g*x^2)^(p + 1), x)}}\\
{\leu0{ - 2*b*(a + b*x)^(m - 1)*(c + d*x)^(n + 1)*(e + f*x + g*x^2)^(p + 1)}}\\
{\leu0{ = 0}},
\bds\bfl
where {\ttt{u = m + 2*p + 1}}.
\efl\ads
{\len0{(m + 2*p + 1)*u*rc}}\\
{\leu1{*INT((a + b*x)^m *(c + d*x)^n *(e + f*x + g*x^2)^p, x)}}\\
{\leu0{ - 2*n*v*d}}\\
{\leu1{*INT((a + b*x)^m *(c + d*x)^(n - 1)*(e + f*x + g*x^2)^(p + 1), x)}}\\
{\leu0{ - 2*(u*b*c - v*a - n*b*d*x)}}\\
{\leu1{*(a + b*x)^(m - 1)*(c + d*x)^n *(e + f*x + g*x^2)^(p + 1) = 0}},
\bds\bfl
where {\ttt{u = m + 2*p + 2}}, {\ttt{v = (m + n + 2*p + 2)*d}}.
\efl

\tttdef\formxii{(a + b*x)^m *(c + d*x)^n *(e + f*x)^p *(g + h*x)^q}
\section{Integrands involving \formxii}

Abbreviations used with these integrands:
\ads
{\reu0{ra = (e*h + f*g)*b - a*f*h}}, {\ttt{rb = b*e*g}},
{\ttt{re = (a*h + b*g)*f - b*e*h}}, {\ttt{rf = a*f*g}}.\\[1ex]
{\reu0{saa = b*c*f*h + ra*d}}, {\ttt{sab = a*c*f*h + ra*c + rb*d}}, 
{\ttt{sbb = rb*c}}.

\subsection{Nondegenerate integrands}
\label{12A}

\begin{sloppypar}
To exclude integrands with confluent roots, the following recurrences 
should be applied only if the overall discriminant does not vanish: 
{\ttt{(a*d - b*c)*(a*f - b*e)*(a*h - b*g)*(c*f - d*e)*(c*h - d*g)*{\fb}%
(e*h - f*g) \tneq 0}}.
\end{sloppypar}
\ads
{\len0{(m + 1)*(a*d - b*c)*(a*f - b*e)*(a*h - b*g)*INT((A + B*x + C*x^2)}}\\
{\leu1{*(a + b*x)^m *(c + d*x)^n *(e + f*x)^p *(g + h*x)^q, x)}}\\
{\leu0{ + INT(((m + 1)*(A*(saa*a - sab*b) + (B*b - C*a)*sbb)}}\\
{\leu1{ - W*((n + 1)*d*e*g + (p + 1)*c*f*g + (q + 1)*c*e*h)}}\\
{\leu1{ - ((m + 1)*((A*b - B*a)*saa + C*(sab*a - sbb*b))}}\\
{\leu1{ + W*((n + p + 2)*d*f*g + (n + q + 2)*d*e*h + (p + q + 2)*c*f*h))*x}}\\
{\leu1{ - (m + n + p + q + 4)*W*d*f*h*x^2)}}\\
{\leu1{*(a + b*x)^(m + 1)*(c + d*x)^n *(e + f*x)^p *(g + h*x)^q, x)}}\\
{\leu0{ + W*(a + b*x)^(m + 1)*(c + d*x)^(n + 1)*(e + f*x)^(p + 1)}}\\
{\leu1{*(g + h*x)^(q + 1) = 0}},
\bds\bfl
where {\ttt{W = A*b^2 - B*a*b + C*a^2}}.
\efl\ads
{\len0{u*INT((A + B*x + C*x^2)*(a + b*x)^m}}\\
{\leu1{*(c + d*x)^n *(e + f*x)^p *(g + h*x)^q, x)}}\\
{\leu0{ - INT(((A*u - V)*a - m*C*sbb + ((A*u - V)*b}}\\
{\leu1{ + (B*u - W)*a - m*C*sab)*x + ((B*u - W)*b - m*C*saa)*x^2)}}\\
{\leu1{*(a + b*x)^(m - 1)*(c + d*x)^n *(e + f*x)^p *(g + h*x)^q, x)}}\\
{\leu0{ - C*(a + b*x)^m *(c + d*x)^(n + 1)*(e + f*x)^(p + 1)*(g + h*x)^(q + 1)}}\\
{\leu0{ = 0}},
\bds\bfl
where {\ttt{V = C*((n + 1)*d*e*g + (p + 1)*c*f*g + (q + 1)*c*e*h)}}, 
{\ttt{W = C*((n + p + 2)*d*f*g + (n + q + 2)*d*e*h + %
(p + q + 2)*c*f*h)}}, {\ttt{u = (m + n + p + q + 3)*d*f*h}}.
\efl\ads
{\len0{u*b*INT((A + B*x + C*x^2)}}\\
{\leu1{*(a + b*x)^m *(c + d*x)^n *(e + f*x)^p *(g + h*x)^q, x)}}\\
{\leu0{ + INT((U*c + V*(n*d*e*g + (p + 1)*c*f*g + (q + 1)*c*e*h) + (U*d}}\\
{\leu1{ - C*u*c + V*(v*c*f*h + (n + p + 1)*d*f*g + (n + q + 1)*d*e*h))*x}}\\
{\leu1{ - (C*u - (v + n)*V*f*h)*d*x^2)}}\\
{\leu1{*(a + b*x)^(m + 1)*(c + d*x)^(n - 1)*(e + f*x)^p *(g + h*x)^q, x)}}\\
{\leu0{ - V*(a + b*x)^(m + 1)*(c + d*x)^n *(e + f*x)^(p + 1)*(g + h*x)^(q + 1)}}\\
{\leu0{ = 0}},
\bds\bfl
where {\ttt{U = (m + 1)*(A*ra*b - (B*b - C*a)*rb)}}, 
{\ttt{V = A*b^2 - (B*b - C*a)*a}}, 
{\ttt{u = (m + 1)*(a*f - b*e)*(a*h - b*g)}}, {\ttt{v = m + p + q + 3}}.
\efl

\subsection{Singly degenerate integrands}
\label{12B}

Here the roots of the first two linear factors are required to coincide 
whence the corresponding resultant must vanish: {\ttt{a*d - b*c = 0}}. To 
exclude additional degeneracies, the recurrences should be applied only 
if also {\ttt{(a*f - b*e)*(a*h - b*g)*(e*h - f*g) \tneq 0}}.
\ads
{\len0{u*(a*f - b*e)*(a*h - b*g)*INT((A + B*x)}}\\
{\leu1{*(a + b*x)^m *(c + d*x)^n *(e + f*x)^p *(g + h*x)^q, x)}}\\
{\leu0{ + INT((u*(A*ra - B*rb) + V*((q + 1)*e*h + (p + 1)*f*g)}}\\
{\leu1{ + (u + p + q + 2)*V*f*h*x)}}\\
{\leu1{*(a + b*x)^(m + 1)*(c + d*x)^n *(e + f*x)^p *(g + h*x)^q, x)}}\\
{\leu0{ - V*(a + b*x)^(m + 1)*(c + d*x)^n *(e + f*x)^(p + 1)*(g + h*x)^(q + 1)}}\\
{\leu0{ = 0}},
\bds\bfl
where {\ttt{V = A*b - B*a}}, {\ttt{u = m + n + 1}}.
\efl\ads
{\len0{(m + n + p + q + 2)*f*h*INT((A + B*x)}}\\
{\leu1{*(a + b*x)^m *(c + d*x)^n *(e + f*x)^p *(g + h*x)^q, x)}}\\
{\leu0{ - INT((U*a - V*rb + (U*b - V*ra)*x)}}\\
{\leu1{*(a + b*x)^(m - 1)*(c + d*x)^n *(e + f*x)^p *(g + h*x)^q, x)}}\\
{\leu0{ - B*(a + b*x)^m *(c + d*x)^n *(e + f*x)^(p + 1)*(g + h*x)^(q + 1) = 0}},
\bds\bfl
where {\ttt{U = (m + n + p + q + 2)*A*f*h - B*((p + 1)*f*g + %
(q + 1)*e*h)}}, {\ttt{V = (m + n)*B}}.
\efl\ads
{\len0{(p + 1)*(a*f - b*e)*(e*h - f*g)*INT((A + B*x)}}\\
{\leu1{*(a + b*x)^m *(c + d*x)^n *(e + f*x)^p *(g + h*x)^q, x)}}\\
{\leu0{ - INT(((p + 1)*(A*re - B*rf) + V*((m + n + 1)*b*g + (q + 1)*a*h)}}\\
{\leu1{ + (m + n + p + q + 3)*V*b*h*x)}}\\
{\leu1{*(a + b*x)^m *(c + d*x)^n *(e + f*x)^(p + 1)*(g + h*x)^q, x)}}\\
{\leu0{ + V*(a + b*x)^(m + 1)*(c + d*x)^n *(e + f*x)^(p + 1)*(g + h*x)^(q + 1)}}\\
{\leu0{ = 0}},
\bds\bfl
where {\ttt{V = A*f - B*e}}.
\efl\ads
{\len0{(m + n + p + q + 2)*b*h*INT((A + B*x)}}\\
{\leu1{*(a + b*x)^m *(c + d*x)^n *(e + f*x)^p *(g + h*x)^q, x)}}\\
{\leu0{ - INT((U*e - p*B*rf + (U*f - p*B*re)*x)}}\\
{\leu1{*(a + b*x)^m *(c + d*x)^n *(e + f*x)^(p - 1)*(g + h*x)^q, x)}}\\
{\leu0{ - B*(a + b*x)^(m + 1)*(c + d*x)^n *(e + f*x)^p *(g + h*x)^(q + 1) = 0}},
\bds\bfl
where {\ttt{U = A*(m + n + p + q + 2)*b*h - B*((m + n + 1)*b*g + %
(q + 1)*a*h)}}.
\efl\ads
{\len0{b*INT((A + B*x)}}\\
{\leu1{*(a + b*x)^m *(c + d*x)^n *(e + f*x)^p *(g + h*x)^q, x)}}\\
{\leu0{ - d*INT((A + B*x)*(a + b*x)^(m + 1)*(c + d*x)^(n - 1)*(e + f*x)^p}}\\
{\leu1{*(g + h*x)^q, x) = 0}}.
\ads
{\len0{u*b*INT((A + B*x)}}\\
{\leu1{*(a + b*x)^m *(c + d*x)^n *(e + f*x)^p *(g + h*x)^q, x)}}\\
{\leu0{ - INT((B*u*e + V*(v*e*h + p*f*g) + (B*u + (v + p)*V*h)*f*x)}}\\
{\leu1{*(a + b*x)^(m + 1)*(c + d*x)^n *(e + f*x)^(p - 1)*(g + h*x)^q, x)}}\\
{\leu0{ + V*(a + b*x)^(m + 1)*(c + d*x)^n *(e + f*x)^p *(g + h*x)^(q + 1) = 0}},
\bds\bfl
where {\ttt{V = A*b - B*a}}, {\ttt{u = (m + n + 1)*(a*h - b*g)}}, 
{\ttt{v = m + n + q + 2}}.
\efl\ads
{\len0{u*f*INT((A + B*x)}}\\
{\leu1{*(a + b*x)^m *(c + d*x)^n *(e + f*x)^p *(g + h*x)^q, x)}}\\
{\leu0{ - INT((B*u*a + V*(v*a*h + (m + n)*b*g) + (B*u + (v + m + n)*V*h)*b*x)}}\\
{\leu1{*(a + b*x)^(m - 1)*(c + d*x)^n *(e + f*x)^(p + 1)*(g + h*x)^q, x)}}\\
{\leu0{ + V*(a + b*x)^m *(c + d*x)^n *(e + f*x)^(p + 1)*(g + h*x)^(q + 1) = 0}},
\bds\bfl
where {\ttt{V = A*f - B*e}}, {\ttt{u = (p + 1)*(e*h - f*g)}}, 
{\ttt{v = p + q + 2}}.
\efl\ads
{\len0{u*f*INT((A + B*x)}}\\
{\leu1{*(a + b*x)^m *(c + d*x)^n *(e + f*x)^p *(g + h*x)^q, x)}}\\
{\leu0{ - INT((B*u*g - V*(v*b*g + q*a*h) + (B*u - (v + q)*V*b)*h*x)}}\\
{\leu1{*(a + b*x)^m *(c + d*x)^n *(e + f*x)^(p + 1)*(g + h*x)^(q - 1), x)}}\\
{\leu0{ - V*(a + b*x)^(m + 1)*(c + d*x)^n *(e + f*x)^(p + 1)*(g + h*x)^q = 0}},
\bds\bfl
where {\ttt{V = A*f - B*e}}, {\ttt{u = (p + 1)*(a*f - b*e)}}, 
{\ttt{v = m + n + p + 2}}.
\efl

\subsection{Doubly degenerate integrands}
\label{12C}

Now the roots of the first three linear factors are required to coincide 
whence the following must vanish: {\ttt{a*d - b*c ={\fb} a*f - b*e = 0}}. 
To exclude additional degeneracies, the recurrences should be applied 
only if also {\ttt{a*h - b*g \tneq 0}}. Note that each relation can be 
used backwards for its own inverse.
\ads
{\len0{u*(a*h - b*g)}}\\
{\leu1{*INT((a + b*x)^m *(c + d*x)^n *(e + f*x)^p *(g + h*x)^q, x)}}\\
{\leu0{ - (u + q + 1)*h}}\\
{\leu1{*INT((a + b*x)^(m + 1)*(c + d*x)^n *(e + f*x)^p *(g + h*x)^q, x)}}\\
{\leu0{ + (a + b*x)^(m + 1)*(c + d*x)^n *(e + f*x)^p *(g + h*x)^(q + 1) = 0}},
\bds\bfl
where {\ttt{u = m + n + p + 1}}.
\efl\ads
{\len0{(q + 1)*(a*h - b*g)}}\\
{\leu1{*INT((a + b*x)^m *(c + d*x)^n *(e + f*x)^p *(g + h*x)^q, x)}}\\
{\leu0{ + (m + n + p + q + 2)*b}}\\
{\leu1{*INT((a + b*x)^m *(c + d*x)^n *(e + f*x)^p *(g + h*x)^(q + 1), x)}}\\
{\leu0{ - (a + b*x)^(m + 1)*(c + d*x)^n *(e + f*x)^p *(g + h*x)^(q + 1) = 0}}.
\ads
{\len0{b*INT((a + b*x)^m *(c + d*x)^n *(e + f*x)^p *(g + h*x)^q, x)}}\\
{\leu0{ - d*INT((a + b*x)^(m + 1)*(c + d*x)^(n - 1)*(e + f*x)^p}}\\
{\leu1{*(g + h*x)^q, x) = 0}}.
\ads
{\len0{(m + n + p + 1)*b*}}\\
{\leu1{INT((a + b*x)^m *(c + d*x)^n *(e + f*x)^p *(g + h*x)^q, x)}}\\
{\leu0{ + q*h*INT((a + b*x)^(m + 1)*(c + d*x)^n *(e + f*x)^p}}\\
{\leu1{*(g + h*x)^(q - 1), x)}}\\
{\leu0{ - (a + b*x)^(m + 1)*(c + d*x)^n *(e + f*x)^p *(g + h*x)^q = 0}}.
\bds

Now the roots of the first two linear factors are required to coincide as 
are the roots of the last two linear factors whence the following must 
vanish: {\ttt{a*d - b*c ={\fb} e*h - f*g = 0}}. To exclude additional 
degeneracies, the recurrences should be applied only if also 
{\ttt{a*f - b*e \tneq 0}}. Note that each relation can be used backwards 
for its own inverse.
\ads
{\len0{(m + n + 1)*(a*f - b*e)}}\\
{\leu1{*INT((a + b*x)^m *(c + d*x)^n *(e + f*x)^p *(g + h*x)^q, x)}}\\
{\leu0{ - (m + n + p + q + 2)*f}}\\
{\leu1{*INT((a + b*x)^(m + 1)*(c + d*x)^n *(e + f*x)^p *(g + h*x)^q, x)}}\\
{\leu0{ + (a + b*x)^(m + 1)*(c + d*x)^n *(e + f*x)^(p + 1)*(g + h*x)^q = 0}}.
\ads
{\len0{b*INT((a + b*x)^m *(c + d*x)^n *(e + f*x)^p *(g + h*x)^q, x)}}\\
{\leu0{ - d*INT((a + b*x)^(m + 1)*(c + d*x)^(n - 1)*(e + f*x)^p}}\\
{\leu1{*(g + h*x)^q, x) = 0}}.
\ads
{\len0{(m + n + 1)*b*}}\\
{\leu1{INT((a + b*x)^m *(c + d*x)^n *(e + f*x)^p *(g + h*x)^q, x)}}\\
{\leu0{ + (p + q)*f*INT((a + b*x)^(m + 1)*(c + d*x)^n *(e + f*x)^(p - 1)}}\\
{\leu1{*(g + h*x)^q, x)}}\\
{\leu0{ - (a + b*x)^(m + 1)*(c + d*x)^n *(e + f*x)^p *(g + h*x)^q = 0}}.

\endgroup 
\section*{Acknowledgments}

The computational advantage and systematic feasibility of the two-term 
recurrences occured to the author in the course of transoceanic 
discussions with Albert D. Rich about the development of his rule-based 
integrator {\sc{Rubi}}. Corresponding formulae for trigonometric 
integrands from the book by Timofeyev and from a table of indefinite 
integrals by W. Gr\"obner {\it{et al.}}\ (Integraltafel, 1.~Teil: 
Unbestimmte Integrale, 1944), reproduced in Gradshteyn{\hd}Ryzhik as 
items 2.551.1 and 2.554.1 and as item 2.558.1, provided seminal 
inspiration. {\sc{Maple}} and {\sc{Mathematica}} results for the example 
integral were supplied by denizens of the {\sss{sci.math.symbolic\ic}} 
newsgroup.

\bibliographystyle{amsplain}

\end{document}